\documentclass{amsart}
\usepackage[utf8]{inputenc}
\usepackage[T1]{fontenc}
\usepackage[colorinlistoftodos]{todonotes}
\usepackage{xifthen, xfrac, mathrsfs, url, enumitem, amssymb, subcaption} 
\usepackage{tikz,pgfplots}
\usepackage{hyperref}
\newtheorem{thm}{Theorem}[section]
\newtheorem{lem}[thm]{Lemma}
\newtheorem{obs}[thm]{Observation}
\newtheorem{prop}[thm]{Proposition}
\newtheorem{alg}{Algorithm}
\theoremstyle{remark}
\newtheorem{rem}{Remark}
\newenvironment{poc}[1][]{\begin{proof}[\ifthenelse{\equal{#1}{}}{Proof of correctness}{Proof of correctness of #1}]}{\end{proof}}

\newcommand{\FF}{\mathbb{F}}
\newcommand{\QQ}{\mathbb{Q}}
\newcommand{\RR}{\mathbb{R}}
\newcommand{\ZZ}{\mathbb{Z}}

\newcommand{\SB}{\mathscr{B}}

\newcommand{\gp}{\mathfrak{p}}
\newcommand{\gpi}{{\mathfrak{p}_i}}
\newcommand{\gq}{\mathfrak{q}}
\newcommand{\gr}{\mathfrak{r}}
\newcommand{\gri}{{\mathfrak{r}_i}}

\newcommand{\GP}{\mathfrak{P}}
\newcommand{\GQ}{\mathfrak{Q}}
\newcommand{\Kp}{K_\gp}

\newcommand{\Kq}{K_\gq}
\newcommand{\Kr}{K_\gr}
\newcommand{\Kri}{K_{\gri}}
\newcommand{\q}{q}
\newcommand{\hq}[1]{\hat{q}_{#1}}
\newcommand{\tq}[1]{\tilde{q}_{#1}}

\newcommand{\st}{\mathrel{\mid}}
\newcommand{\abs}[2][\gr]{|#2|_{#1}}

\newcommand{\form}[1]{\langle#1\rangle}
\newcommand{\ext}[2]{#1/#2}
\newcommand{\quo}[2]{\sfrac{#1}{#2}}
\newcommand{\Sing}[1]{\mathbb{E}_{#1}}

\newcommand{\SingS}{\Sing{S}}

\newcommand{\SingSq}{\Sing{S\cup \{\gq\}}}

\newcommand{\Units}[1]{\mathbb{U}_{#1}}
\newcommand{\units}[1]{#1^\times}
\newcommand{\un}[1][]{\ifthenelse{\equal{#1}{}}{\units{K}}{\units{#1}}}
\newcommand{\squares}[1]{#1^{\times2}}
\newcommand{\sq}[1][]{\ifthenelse{\equal{#1}{}}{\squares{K}}{\squares{#1}}}
\newcommand{\sqgd}[1][]{\ifthenelse{\equal{#1}{}}{\sfrac{\un}{\sq}}{\sfrac{\units{#1}}{\squares{#1}}}}
\newcommand{\norm}[2][\ext{L}{K}]{\ifthenelse{\equal{#2}{}}{N_{#1}}{N_{#1}\left(#2\right)}}

\DeclareMathOperator{\disc}{disc}
\DeclareMathOperator{\ord}{ord}
\DeclareMathOperator{\sgn}{sgn}

\newcommand{\class}[2][]{\ifthenelse{\equal{#1}{}}{[#2]}{[#2]_{#1}}}   

\newcommand{\into}{\rightarrowtail}

\newcommand{\term}[1]{\emph{#1}}
\makeatletter
\@namedef{subjclassname@2020}{\textup{2020} Mathematics Subject Classification}
\makeatother
\author{Przemysław Koprowski}
\email{przemyslaw.koprowski@us.edu.pl}
\address{Institute of Mathematics, University of Silesia, ul. Bankowa 14, Katowice, Poland}
\title{Isotropic vectors over global fields}
\begin{document}
\begin{abstract}
An isotropic vector of a given quadratic form is a nonzero vector where the form vanishes. Geometrically, it is a vector that is self-orthogonal with respect to this form. On the other hand, from the arithmetical point of view, an isotropic vector forms a solution to a multivariate quadratic equation. The problem of constructing isotropic vectors is one of the leading forces in the computational theory of quadratic forms. In this paper, we present algorithms for finding isotropic vectors of quadratic forms (of any dimension) over an arbitrary global field of characteristic distinct from~$2$.
\end{abstract}
\maketitle

\section{Introduction} 
The notion of isotropy is fundamental for the algebraic theory of quadratic forms. Given a quadratic form~$\q$, a nonzero vector~$v$ is called \term{isotropic} when $\q(v) = 0$. The problem of finding an isotropic vector of a quadratic form has quite a long history with several ramifications. The first effective method (in fact, an algorithm, although it was not called as such when it was invented) is a descent method for ternary quadratic forms over the rationals. It dates back to the 18th century. The reader can find its description in an explicit algorithmic form in, e.g., \cite[Chapter~IV]{Smart1998}, where it is attributed to Lagrange, or in \cite[Algorithm~I]{CR2003}, where it is attributed to Legendre. For a long time, it was the standard method for finding isotropic vectors of ternary quadratic forms over~$\QQ$. The advent of computational algebra in the second half of the last century brought the problem back to the spotlight and stipulated the development of new, more efficient algorithms. It is well known (and recalled in Section~\ref{sec:dim3}) that finding a solution of a ternary form is equivalent to solving a norm equation in a quadratic field extension. The latter problem can be considered for arbitrary (i.e., not only quadratic) field extensions and has been a subject of separate studies. We shall refer the reader to \cite{Cohen2000, FJP1997, FP1983, Garbanati1980, LJS2024, Simon2002} for some known results. In particular, it should be pointed out that Simon's solution based on $S$-units (see \cite{Simon2002}) served as an original inspiration to our algorithm for finding isotropic vectors of quaternary quadratic forms (see Algorithm~\ref{alg:dim4}).

The main disadvantage of the Lagrange/Legendre method is that it requires multiple factorizations of integers. While for low-dimensional forms (i.e., ternary and quaternary forms), it is necessary to factor the discriminant of the form, no other factorizations are necessary. In 1998, Cochrane and Mitchel \cite{CM1998} and in 2003, Cremona and Rusin \cite{CR2003} proposed algorithms for ternary quadratic forms over~$\QQ$ that need no other factorization except the factorization of the discriminant. For forms over~$\QQ$ that are not diagonal, Simon \cite{Simon2005} devised an algorithm that finds an isotropic vector without first diagonalizing the form and also using only the factorization of the discriminant. For forms over the rationals of dimensions greater than~$4$, it is possible altogether to avoid \emph{any} factorization of integers as it was shown by Castel \cite{Castel2013}. The key idea of Castel's solution is first to find an indefinite subform of dimension~$5$ and then extend it to a $6$-dimensional isotropic form of a prescribed discriminant.

Far less is known for forms over fields different from~$\QQ$. Schicho \cite{Schicho1998, Schicho2000} generalized the descent method to quadratic forms over the field $\RR(t)$ of rational functions over the reals. The algorithm was further generalized to an arbitrary rational function field $K(t)$ of characteristic distinct from~$2$ by van Hoeij and Cremona in \cite{vHC2006}. A different algorithm that works with quadratic forms over $\FF_q(t)$ for odd~$q$ was given by Ivanyos et al. in \cite{IKR2019}. Recently, Kutas et al. \cite{KMZC2022} presented a method for finding an isotropic vector of a quaternary quadratic form over rational function fields in characteristic~$2$. Finally, we should also mention that it is possible to generalize the descent method to ternary quadratic forms over number fields. A corresponding algorithm is implicit in \cite{Pohst2000} and was implemented in Magma (see \cite[\S127.5]{CBFS2023}) by Steve Donelly. To the best of our knowledge, it has never been published in an explicit algorithmic form in any research paper. 

In general, the methods listed above that work with forms over the rationals heavily depend on the idiosyncrasies of the arithmetic of the (euclidean) ring of integers. Consequently, they do not directly generalize to number fields. Here, the problem of constructing an isotropic vector is equivalent to solving a degree~$2$ homogeneous multivariate Diophantine equation over a ring of algebraic integers, where the arithmetic of such ring can be quite\dots{} weird. That is probably the reason why the problem of finding isotropic vectors of general quadratic forms over number fields has---to the best of our knowledge---remained mostly unsolved. The sole purpose of the present paper is to remedy this situation to a certain degree. We present explicit algorithms that construct isotopic vectors for quadratic forms over arbitrary global fields of characteristic distinct from~$2$. The algorithms described in this paper have been implemented in Magma \cite{BCP1997} and can be downloaded from the author's homepage at \url{http://www.pkoprowski.eu/cqf}. In Appendix~\ref{sec:implementation}, we report how they perform in practice.

\section{Notation}
Throughout this paper, we use the following notation and terminology. Since the theory of quadratic forms in characteristic~$2$ is very different from all other characteristics, we exclude the case of characteristic~$2$ from all our considerations. For this reason, in the rest of the paper, the term ``global field'' actually means ``global field of characteristic distinct from~$2$''. By~$K$, we will denote our base field, which is either a number field or an algebraic function field over some finite field~$\FF_q$, where $q$ is a power of an odd prime. If $\ext{L}{K}$ is a finite field extension, then $\norm{}: L \to K$ is the associated norm. A diagonal quadratic form $a_1x_1^2 + \dotsb + a_nx_n^2$ is denoted $\form{a_1, \dotsc, a_n}$. Since the problem of constructing isotropic vectors of degenerate forms is trivial and so totally uninteresting, all forms are always assumed to be non-degenerate, meaning that every $a_i$ is nonzero.

A \term{place} of~$K$ is an equivalence class of valuations. We denote the places using fraktur letters $\gp$, $\gq$, $\gr$. A place of a number field can be either Archimedean or non-archimedean. Over global function fields, all places are non-archimedean. We write $\Kp$ for the completion of~$K$ at~$\gp$ and $K(\gp)$ for its residue field. In addition, if $\gp$ is non-archimedean, we write $\ord_\gp: \un\to \ZZ$ for the associated normalized discrete valuation. It is well known that instead of valuations, one can alternatively use \term{norms}. Given a place~$\gp$ of~$K$, we will write $\abs[\gp]{\cdot}: K\to \RR_{\geq 0}$ for the associated norm. In fact, in this paper, we will use norms as more convenient objects exclusively in the case of real places. If $K$ is a number field, $\gr$ its real place, and $\sigma_\gr: K\into \RR$ the associated real embedding, then for any $a\in K$ its norm $\abs{a}$ is nothing else but the standard absolute value $|\sigma_\gr(a)|$ of the image of~$a$ under $\sigma_\gr$. 

The non-archimedean places correspond to prime ideals in a maximal order of~$K$. For this reason, we will use the terms ``non-archimedean places'' and ``primes'' interchangeably. We will sometimes write ``finite primes'' to emphasize the fact that they are non-archimedean. A prime~$\gp$ of a number field is called \term{dyadic} if $\ord_\gp 2\neq 0$. Otherwise, it is called \term{non-dyadic}. In global function fields, all primes are non-dyadic. If $\gp$ is a finite prime and $a, b\in \un$, then $(a,b)_\gp$ denotes the \term{Hilbert symbol} of~$a$ and~$b$ at~$\gp$. If further $\q = \form{a_1, \dotsc , a_n}$ is a quadratic form, then following \cite{Lam2005} we define the \term{Hasse invariant} of~$\q$ at~$\gp$ by the formula
\[
s_\gp(\q) := \prod_{i<j} (a_i, a_j)_\gp.
\]

Recall that a finite set~$S$ of places of a global field is called a \term{Hasse set} if it contains all Archimedean places. In particular, every finite set of places of a global function field is a Hasse set. With a quadratic form $\q = \form{a_1, \dotsc , a_n}$ we associate a Hasse set, denoted~$\GP(\q)$, consisting of:
\begin{itemize}
\item all Archimedean places and all dyadic primes (provided that $K$ is a number field),
\item all these primes $\gp$ such that $\ord_\gp a_i$ is odd for at least one $i\leq n$.
\end{itemize}
Let $S$ be a Hasse set, an element $a\in \un$ is said to be \term{$S$-singular} if $\ord_\gp a$ is even for every prime~$\gp$, not in~$S$. It is called an \term{$S$-unit} if $\ord_\gp a = 0$ for every $\gp\notin S$. The group of all $S$-singular elements is denoted $E_S$, and the group of all $S$-units by $U_S$. Observe that $\sq$ is a subgroup of $E_S$, hence we may consider the quotient group $\SingS := \quo{E_S}{\sq}$ of $S$-singular square classes. It is a subgroup of the square class group $\sqgd$ of~$K$. Moreover, we will denote $\Units{S} := \quo{U_S}{U_S\cap \sq}$ the group of $S$-units modulo squares. It can be treated as a subgroup of $\SingS$. Observe that $\SingS$ is an elementary $2$-group. Hence, it is an $\FF_2$-linear space.

\section{Algorithmic prerequisites}\label{sec:prereqs}
The algorithms presented in this paper depend on the number of already existing tools. First of all, we need routines for solving systems of linear equations (basically over $\FF_2$) and for basic arithmetic in a global field (see, e.g., \cite{Cohen1993}). As far as more advanced tools are concerned, we need:
{
\setlist[enumerate]{wide=\parindent}
\begin{enumerate}
\item Factorization of ideals in global fields (see, e.g., \cite{Cohen1993, GMN2013}). 
\item Solution to a norm equation of a form $\norm{x} = b$, where $b\in \un$ and $\ext{L}{K}$ is a finite field extension. As was mentioned in the introduction, this problem is discussed in \cite{Cohen2000, FJP1997, FP1983, Garbanati1980, LJS2024, Simon2002}.
\item Construction of a basis (over $\FF_2$) of the group $\SingS$ of $S$-singular square classes. The author is aware of three existing methods. One is described in \cite{Koprowski2021}. Another one is outlined in the Magma manual \cite{CBFS2023}. In this method, one constructs a set~$T$ containing~$S$, such that the $T$-class number is odd. Then, it is known that the groups $\Sing{T}$ and $\Units{T}$ coincide. Hence, the sought group~$\SingS$ can be constructed as a subspace of~$\Units{T}$. A third method, due to A.~Czogała, was recently described in \cite{KR2023}.
\end{enumerate}}

\section{Algorithms}
In this section, we present the main results of the paper, namely the algorithms for constructing an isotropic vector of a given quadratic form. Our general strategy is as follows. For forms of dimensions~$2$ or~$3$, the corresponding vector is constructed directly. On the other hand, for higher dimensional forms, we will reduce the problem to the one for forms of lower dimension. Hence, we will recursively push it down the ladder till we reach ternary forms. Unfortunately, the strategy differs significantly from one dimension to the other. Consequently, we present dedicated algorithms for forms of dimensions~$4$ and~$5$ and a general algorithm that deals with forms of dimension~$6$ and above.

\subsection{Dimension 3 and below}\label{sec:dim3}
First, let us deal with forms of dimensions~$3$ and below. Unary forms cannot be isotropic, and so we ignore them. The construction of an isotropic vector of a binary form is trivial and boils down to computing the square root of its discriminant. As explained earlier, for ternary forms, the problem is equivalent to the well-studied problem of solving the norm equation. Indeed, let $\q =\form{a_1, a_2, a_3}$ be an isotropic ternary form and $v = (v_1, v_2, v_3) \in K^3$ be its isotropic vector, we are looking for. Without loss of generality, we may assume that all three binary subforms $\form{a_i, a_j}$, with $1\leq i< j\leq 3$ are anisotropic. In particular, this means that $-\sfrac{a_2}{a_1}$ is not a square in~$K$ and that $v_3$ is non-zero. We have $a_1v_1^2 + a_2v_2^2 + a_3v_3^2 = 0$. Rearranging the terms we write
\[
-\frac{a_3}{a_1} 
= \left(\frac{v_1}{v_3}\right)^2 + \frac{a_2}{a_1}\cdot \left(\frac{v_2}{v_3}\right)^2.
\]
Now, the right-hand side is the norm of an element 
\[
\xi = \frac{v_1}{v_3} + \frac{v_2}{v_3}\cdot \sqrt{-\frac{a_2}{a_1}}
\]
in the quadratic field extension $\ext{K\bigl(\sqrt{-\sfrac{a_2}{a_1}}\bigr)}{K}$. Therefore, as said, the problem of constructing an isotropic vector of a ternary form is equivalent to the one of solving the norm equation.

\subsection{Dimension 4}\label{sec:dim4}
We shall now deal with quaternary forms. This case requires some preparations. Fix an isotropic quaternary quadratic form $\q = \form{a_1, a_2, a_3, a_4}$ over~$K$. Without loss of generality we can assume that no $3$-dimensional subform $\form{a_i, a_j, a_k}$ of~$\q$ (for some $1\leq i\leq j\leq l\leq 4$) is isotropic. Split the form into two subforms $\q_1 := \form{a_1, a_2}$ and $\q_2 := \form{a_3, a_4}$. Let us begin with the following elementary observation determining our strategy for constructing an isotropic vector.

\begin{obs}\label{obs:quaternary_iso_condition}
The form~$\q$ is isotropic if and only if there is an element $b\in\un$ such that the forms $\tq1 := \form{-b}\perp \q_1$ and $\tq2 := \form{b}\perp \q_2$ are both isotropic. Moreover, if $u = (u_1, u_2, u_3)$ is an isotropic vector for~$\tq1$ and $w = (w_1, w_2, w_3)$ is an isotropic vector for~$\tq2$, then $v := (\sfrac{u_2}{u_1}, \sfrac{u_3}{u_1}, \sfrac{w_2}{w_1}, \sfrac{w_3}{w_1})$ is an isotropic vector for~$\q$.
\end{obs}

Consequently, our goal is to find a method for constructing an element~$b$ as above. We shall find it among representatives of $S$-singular square classes for a certain Hasse set~$S$. Let $S_\infty = \{\gr_1, \dotsc, \gr_r\}$ be the set (possibly empty) of all these real places~$\gr$ of~$K$ where either $\sgn_\gr a_1 = \sgn_\gr a_2$ (hence $\q_1$ is anisotropic at~$\gr$) or $\sgn_\gr a_3 = \sgn_\gr a_4$ (hence $\q_2$ is anisotropic at~$\gr$).

\begin{lem}\label{lem:quaternary_real}
For every $\gr\in S_\infty$ the forms $\tq1$ and $\tq2$ are isotropic at~$\gr$ if and only~if 
\begin{equation}\tag{$\spadesuit$}\label{eq:quaternary_real_condition}
\sgn_\gr b =
\begin{cases}
\sgn_\gr a_1, & \text{if $\q_1\otimes \Kr$ is anisotropic,}\\
-\sgn_\gr a_3, & \text{otherwise.}
\end{cases}
\end{equation}
\end{lem}

\begin{proof}
First, suppose that $\q_1$ is anisotropic at~$\gr$. For $\tq1\otimes \Kr$ to be isotropic, we must have $\sgn_\gr(-b) \neq \sgn_\gr a_1$. On the other hand, if $\q_1\otimes \Kr$ is isotropic, it is $\q_2\otimes \Kr$ that is anisotropic. Now, for $\tq2\otimes \Kr$ to be isotropic we must have $\sgn_\gr b \neq \sgn_\gr a_3$.

Conversely, suppose that condition~\eqref{eq:quaternary_real_condition} holds. If $\q_1\otimes \Kr$ is anisotropic, then it is clear that $\form{-b}\perp \q_1$ is isotropic at~$\gr$. We must show that $\tq2\otimes \Kr$ is also isotropic. If already the form $\q_2\otimes \Kr$ is isotropic, then so is $\tq2\otimes \Kr$. Hence, suppose that $\q_2\otimes \Kr$y is anisotropic. The form~$\q$ is isotropic by assumption. In particular, it is isotropic at~$\gr$. Therefore, all four coefficients cannot have the same sign. It follows that 
\[
\sgn_\gr b = \sgn_\gr a_1 = \sgn_\gr a_2\neq \sgn_\gr a_3 = \sgn_\gr a_4
\]
and so $\tq2\otimes \Kr$ is isotropic as expected. When $\q_1\otimes \Kr$ is isotropic, we need only to show that $\tq2$ is isotropic at~$\gr$. This fact, however, follows immediately from the second part of condition~\eqref{eq:quaternary_real_condition}.
\end{proof}

The next lemma is an analog of Lemma~\ref{lem:quaternary_real} but for finite primes.

\begin{lem}\label{lem:quaternary_finite}
For every finite prime~$\gp$:
\begin{enumerate}
\item The form~$\tq1$ is locally isotropic at~$\gp$ if and only if $(a_1, a_2)_\gp = (-a_1a_2, b)_\gp$. 
\item The form~$\tq2$ is locally isotropic at~$\gp$ if and only if $(-a_3, -a_4)_\gp = (-a_3a_4, b)_\gp$.
\end{enumerate}
\end{lem}

\begin{proof}
It suffices to prove just one of the two assertions. The other one is fully analogous. Now, \cite[Proposition~V.3.22]{Lam2005} asserts that a ternary form~$\tq1$ is isotropic over a local field~$\Kp$ if and only if its Hasse invariant equals $( -1, -\det \tq1)_\gp$. Therefore, the isotropy of $\tq1\otimes \Kp$ is equivalent to the condition
\[
s_\gp(\tq1) 
= (-b, a_1)_\gp(-b, a_2)_\gp(a_1, a_2)_\gp 
= (-1, a_1a_2b)_\gp.
\]
The Hilbert symbol is symmetric and multiplicative with respect to each of its arguments. Thus, the above condition can be rewritten in the form $(a_1, a_2)_\gp = (-a_1a_2, b)_\gp$.
\end{proof}

It follows from the Local-Global Principle (see, e.g., \cite[Section~VI.3]{Lam2005}) that $\tq1$ and $\tq2$ are isotropic (over~$K$) if and only if the conditions of Lemmas~\ref{lem:quaternary_real} and \ref{lem:quaternary_finite} hold for all places of~$K$. Now, for complex places, real places not in~$S_\infty$ and finite primes that do not divide any of the elements $a_1, \dotsc, a_4, b$, these two conditions are trivially satisfied. Thus, we have only finitely many places to deal with. Fix a Hasse set~$S$ containing $\GP(\q)$ and such that $b$ is $S$-singular. Let $\SB = \{\kappa_1, \dotsc, \kappa_k\}$ be a basis of~$\SingS$ treated as a linear space over~$\FF_2$. Write the square class of~$b$ as a linear combination of elements from~$\SB$: 
\[
b\equiv \kappa_1^{\xi_1}\dotsm \kappa_k^{\xi_k}\pmod{\sq},
\]
for some $\xi_1, \dotsc, \xi_k\in \{0,1\}$. We will now define three vectors $u$, $v$, $w$ and corresponding matrices $A$, $B$, $C$ over~$\FF_2$. To this end, we proceed as follows.

For every real place $\gri\in S_\infty$ let $u_i\in \{0,1\}$ be such that $(-1)^{u_i} = \sgn_\gri a_1$, when $\q_1\otimes \Kri$ is anisotropic, and $(-1)^{u_i} = -\sgn_\gri a_3$, when $\q_1\otimes \Kri$ is isotropic and so $\q_2\otimes \Kri$ is anisotropic. Further, take $\alpha_{ij}\in \{0, 1\}$ such that $(-1)^{\alpha_{ij}} = \sgn_\gri \kappa_j$. Define $A = (\alpha_{ij})$ to be a $r\times k$ matrix with rows indexed by places in~$S_\infty$ and columns by elements of~$\SB$. Finally, let $u = (u_i)$ be a column vector indexed by real places in~$S_\infty$.

Analogously, for every finite prime $\gp_i\in S$, let $v_i, w_i\in \{0, 1\}$ be such that 
\[
(-1)^{v_i} = (a_1, a_2)_\gpi
\qquad\text{and}\qquad
(-1)^{w_i} = (-a_3, -a_4)_\gpi.
\]
Let $\beta_{ij}, \gamma_{ij}\in \{0, 1\}$ be such that 
\[
(-1)^{\beta_{ij}} = (-a_1a_2, \kappa_j)_\gpi
\qquad\text{and}\quad
(-1)^{\gamma_{ij}} = (-a_3a_4, \kappa_j)_\gpi.
\]
Define the matrices $B := (\beta_{ij})$, $C = (\gamma_{ij})$ with rows indexed by finite primes in~$S$ and columns by elements of the basis~$\SB$. Likewise define column vectors $v := (v_i)$ and $w := (w_i)$.

\begin{lem}\label{lem:quaternary_linear_system}
Use the notation introduced above. The coordinates $\xi_1, \dotsc, \xi_k\in \FF_2$ of~$b$ with respect to the basis~$\SB$ of~$\SingS$ form a solution to the following system of $\FF_2$-linear equations
\begin{equation}\tag{$\clubsuit$}\label{eq:quaternary}
\left(\begin{array}{c} A\\\hline B\\\hline C \end{array}\right)
\cdot 
\begin{pmatrix} x_1\\ \vdots\\ x_k\end{pmatrix}
=
\left(\begin{array}{c} u\\\hline v\\\hline w \end{array}\right).
\end{equation}
In particular, an $S$-singular element~$b$ exists if and only if system~\eqref{eq:quaternary} is solvable.
\end{lem}

\begin{proof}
We have $b = \kappa_1^{\xi_1}\dotsm \kappa_k^{\xi_k}\cdot c^2$ for some $c\in \un$. Fix a real place $\gri\in S_\infty$, and assume that $\q_1\otimes \Kri$ is anisotropic. Lemma~\ref{lem:quaternary_real} asserts that $\tq1\otimes \Kri$ and $\tq2\otimes \Kri$ are isotropic if and only if
\begin{multline*}
(-1)^{u_i}
= \sgn_\gri a_1
= \sgn_\gri b
\\
= (\sgn_\gri \kappa_1)^{\xi_i}\dotsm (\sgn_\gri \kappa_k)^{\xi_k}
= (-1)^{\alpha_{i1}\xi_1}\dotsm (-1)^{\alpha_{ik}\xi_k}.
\end{multline*}
This gives us the following $\FF_2$-linear equation:
\[
\sum_{j\leq k} \alpha_{ij}\xi_j = u_i.
\]
In the case when is is $\q_2\otimes \Kri$ that is anisotropic an analogous argument leads to precisely the same linear equation.

Now, fix a prime $\gpi\in S$. Using Lemma~\ref{lem:quaternary_finite} we can write 
\begin{multline*}
(-1)^{v_i}
= (a_1, a_2)_\gpi
= (-a_1a_2, b)_\gpi
\\
= (-a_1a_2, \kappa_1)_\gpi^{\xi_1}\dotsm (-a_1a_2, \kappa_k)_\gpi^{\xi_k}
= (-1)^{\beta_{i1}\xi_1}\dotsm (-1)^{\beta_{ik}\xi_k}.
\end{multline*}
This gives us the next equation:
\[
\sum_{j\leq k} \beta_{ij}\xi_j = v_i.
\]
Analogously, we obtain the last equation:
\[
\sum_{j\leq k} \gamma_{ij}\xi_j = w_i.
\]
Writing all three equations in a matrix form for all possible indices~$i$, we obtain the assertion of the lemma. 
\end{proof}

One possible method for finding an isotropic vector of a quaternary form is now apparent. Start with a Hasse set $S = \GP(\q)$ and enlarge it prime-by-prime till system~\eqref{eq:quaternary} becomes solvable. The next proposition shows that we do not have to add too many primes.

\begin{prop}\label{prop:point_q}
Keep the above notation. There is a prime~$\gq$ of~$K$ such that for $S' = S\cup \{\gq\}$ system~\eqref{eq:quaternary} is solvable.
\end{prop}

\begin{proof}
The form~$\q$ is isotropic by assumption. Therefore, there is some (not necessarily $S'$-singular) element~$b$ satisfying the assertions of  Observation~\ref{obs:quaternary_iso_condition}. Now, \cite[Lemma~2.1]{LW1992} says that there is a prime $\gq\notin S$ and a $\bigl(S\cup \{\gq\}\bigr)$-singular element~$c$ such that
\[
c\equiv b \pmod{\gp^{1 + \ord_\gp 4}} .
\]
for every prime $\gp\in S$ and $\sgn_\gr c = \sgn_\gr b$ for all real places~$\gr$ of~$K$. It follows from the Local Square Theorem (see, e.g., \cite[Theorem VI.2.19]{Lam2005}) that the local square classes $b\cdot \squares{\Kp}$ and $c\cdot \squares{\Kp}$ coincide for every place $\gp\in S$. We claim that $c$ is represented over~$K$ by the forms $\form{a_1, a_2}$ and $\form{-a_3,-a_4}$. 

Denote $\hq1 := \form{-c, a_1, a_2}$ and $\hq2 := \form{c, a_3, a_4}$. Evidently, the forms~$\hq1$ and~$\hq2$ are locally isotropic at all complex places of~$K$. Fix any place $\gp\in S$, be it Archimedean or non-archimedean. We know that the forms $\form{-b, a_1, a_2}$ and $\form{b, a_3, a_4}$ are isotropic, so are their localizations $\form{-b, a_1, a_2}\otimes \Kp$ and $\form{b, a_3, a_4}\otimes \Kp$. Now, the local square classes $b\cdot \squares{\Kp}$ and $c\cdot \squares{\Kp}$ are equal. Hence,
\[
\form{-b, a_1, a_2}\otimes \Kp \cong \hq1\otimes \Kp
\qquad\text{and}\qquad
\form{b, a_3, a_3}\otimes \Kp \cong \hq2\otimes \Kp.
\]
This means that $\hq1\otimes \Kp$ and $\hq2\otimes \Kp$ are both isotropic, as desired.

Next, take a finite prime~$\gp$ not in~$S$ but distinct from~$\gq$. In particular, $\gp$ is non-dyadic since all the dyadic primes are in $\GP(\q)\subseteq S$. Then all the elements $a_1, \dots, a_4$ and~$c$ have even valuations at~$\gp$, consequently the forms $\hq1\otimes \Kp$ and $\hq2\otimes \Kp$ are isotropic by \cite[Corollary~VI.2.5]{Lam2005}. 

Finally, consider the singled-out prime~$\gq$. By the previous paragraph and \cite[Proposition~V.3.22]{Lam2005} for every prime $\gp\neq \gq$ the Hasse invariant $s_\gp(\hq1)$ equals the Hilbert symbol $(-1, -\det \hq1)_\gp$. Thus, we have 
\[
1 
= s_\gp(\hq1)\cdot (-1 -\det \hq1)_\gp
= (c, -a_1a_2)_\gp\cdot (a_1, a_2)_\gp.
\]
Using Hilbert reciprocity law (see e.g., \cite[Chapter~VII]{OMeara2000}) we can write
\begin{align*}
1
&= (c, -a_1a_2)_\gq\cdot (a_1, a_2)_\gq\cdot \prod_{\gp\neq \gq}(c, -a_1a_2)_\gp\cdot \prod_{\gp\neq \gq} (a_1,a_2)_\gp\\
&= (c, -a_1a_2)_\gq\cdot (a_1, a_2)_\gq\\
&= s_\gq(\hq1)\cdot (-1, -\det\hq1)_\gq.
\end{align*}
Hence, again by \cite[Proposition~V.3.22]{Lam2005} we obtain that $\hq1\otimes \Kq$ is isotopic. A similar argument shows that $\hq2\otimes \Kq$ is isotropic, as well. We have shown that the forms~$\hq1$ and~$\hq2$ are isotropic at every completion of~$K$. Therefore, by the Local-Global Principle, they are isotropic over~$K$, as well. This proves the claim. It now follows from the previous lemma that system~\eqref{eq:quaternary} is solvable in~$\SingSq$.
\end{proof}

We are now ready to present the main algorithm of this paper.

\begin{alg}\label{alg:dim4}
Given an isotropic quadratic form $\q = \form{a_1, a_2, a_3, a_4}$ of dimension~$4$ over a global field~$K$, this algorithm outputs a vector $v\in K^4$ such that $\q(v) = 0$.
\begin{enumerate}
\item\label{st:dim4:quick_exit} If there are indices $1\leq i < j < k \leq 4$ such that $\form{a_i, a_j, a_k}$ is isotropic, then construct its isotropic vector $v\in K^3$. Output $v$ with zero inserted at the omitted spot and quit.
\item\label{st:dim4:initilization} Initialize $S := \GP(\q)$, and denote $\q_1 := \form{a_1, a_2}$ and $\q_2 := \form{a_3, a_4}$.
\item\label{st:dim4:real_places} Let $S_\infty := \{\gr_1, \dotsc, \gr_r\}\subset S$ be the set of all the real places of~$K$ at which at least one of the forms $\q_1$, $\q_2$ is anisotropic:
\[
S_\infty := \bigl\{ \gr\st \sgn_\gr a_1 = \sgn_\gr a_2\text{ or }\sgn_\gr a_3 = \sgn_\gr a_4\bigr\}.
\]
\item\label{st:dim4:vec_u} Let $u = (u_i)$ be the \textup(column\textup) vector defined by the formula
\[
(-1)^{u_i} =
\begin{cases}
\sgn_\gri a_1 & \text{if $\q_1\otimes \Kri$ is anisotropic}\\
-\sgn_\gri a_3 &\text{otherwise.}
\end{cases}
\]
\item\label{st:dim4:loop} Repeat the following steps: 
  \begin{enumerate}
  \item\label{st:dim4:basis} Construct a basis $\SB = \{ \kappa_1, \dotsc, \kappa_k\}$ of~$\SingS$. 
  \item\label{st:dim4:vec_v} Construct a \textup(column\textup) vector $v = (v_i)$ determined by the condition
  \[
  (-1)^{v_i} = (a_1, a_2)_\gpi.
  \]
  \item\label{st:dim4:vec_w} Construct a \textup(column\textup) vector $w = (w_i)$ given by the formula
  \[
  (-1)^{w_i} = (-a_3, -a_4)_\gpi.
  \]
  \item\label{st:dim4:mtx_A} Construct a matrix $A = (\alpha_{ij})$ defined by the condition
  \[
  (-1)^{\alpha_{ij}} = \sgn_\gri \kappa_j.
  \]
  \item\label{st:dim4:mtx_B} Construct a matrix $B = (\beta_{ij})$ determined by the formula
  \[
  (-1)^{\beta_{ij}} = ( -a_1a_2, \kappa_j)_\gpi.
  \]
  \item\label{st:dim4:mtx_C} Construct a matrix $C = (\gamma_{ij})$ given by the condition
  \[
  (-1)^{\gamma_{ij}} = ( -a_3a_4, \kappa_j)_\gpi.
  \]
  \item Form system~\eqref{eq:quaternary}, by concatenating matrices $A$, $B$, $C$ and vectors $u$, $v$, $w$. Check if it is solvable.
  \item If it is, denote a solution by $(\xi_1, \dotsc, \xi_k)$ and exit the loop. 
  \item\label{st:dim4:new_prime} Otherwise, if system~\eqref{eq:quaternary} has no solution, find a new prime $\gq\notin S$. Append it to~$S$ and reiterate the loop.
  \end{enumerate}
\item Set 
\[
b := \kappa_1^{\xi_1}\dotsm \kappa_k^{\xi_k},\qquad
\hq1 := \form{-b, a_1, a_2},\qquad
\hq2 := \form{ b, a_3, a_4}.
\]
\item\label{st:dim4:solve_ternary} Find isotropic vectors of ternary forms~$\hq1$ and~$\hq2$. Denote them~$u$ and~$w$, respectively. 
\item Output $v := (\sfrac{u_2}{u_1}, \sfrac{u_3}{u_1}, \sfrac{w_2}{w_1}, \sfrac{w_3}{w_1})$.
\end{enumerate}
\end{alg}

\begin{rem}\label{rem:incremental}
If one uses an incremental algorithm for constructing a basis of a group~$\SingS$ of $S$-singular square classes (as described, e.g., in \cite{Koprowski2021}), then it is not necessary to compute the basis~$\SB$ in step~\eqref{st:dim4:basis} from scratch every time. Instead, one may just update the basis computed in the previous iteration of the loop. 
\end{rem}

\begin{rem}\label{rem:density}
In the proof of Proposition~\ref{prop:point_q}, we used \cite[Lemma~2.1]{LW1992} to prove the existence of a point~$\gq$ such that system~\eqref{eq:quaternary} is solvable in~$\SingSq$. The proof of this lemma relies on the Chebotarev Density Theorem to show that the set of such points is not empty. In particular, this set have a positive density. It means that the corresponding primes~$\gq$ are (in a certain sense) relatively common.
\end{rem}

\begin{rem}\label{rem:exhaustive_search}
In order to rigorously prove that after finitely many iterations, we will come across a prime~$\gq$ satisfying assertions of Proposition~\ref{prop:point_q}, we should exhaustively iterate over all the primes. If~$K$ is a number field, we may, for example, start with a prime number $p = 3$ and loop over all the primes of~$K$ dividing~$p$. Then take $p = 5$ and again loop over all the primes of~$K$ over~$p$, and so on. A similar exhaustive search can be applied to global function fields, as well. In practice, a much more efficient solution is just to pick~$\gq$ at random, using the fact that the set of such primes has a positive density (see the previous remark). Thus, depending on the selected strategy for finding primes in step~\eqref{st:dim4:new_prime} the above algorithm can be viewed either as a deterministic or a randomized one.
\end{rem}

We comment on the complexity of Algorithm~\ref{alg:dim4} in Appendix~\ref{sec:complexity}.

\subsection{Dimension 5 and above}
We shall now deal with forms of dimensions~$5$ and higher. Such forms are necessarily locally isotropic at all finite primes. The general strategy we will use is to reduce the dimension of the form step-by-step till we drop it to~$4$, where we may use Algorithm~\ref{alg:dim4}. Essentially we will use a variant of the $2$-linkage of quadratic forms (see, e.g., \cite[Chapter~8]{Szymiczek1997}) Recall that for any nonzero $a,b,s,t\in \un$ such that $a + b\neq 0$ we have the following isometries of binary quadratic forms (see, e.g., \cite[Corollary~6.4.1]{Szymiczek1997}): 
\[
\form{a, b} \cong \form{as^2, bt^2}
\qquad\text{and}\qquad
\form{a, b} \cong \form{a + b, (a+b)\cdot ab}.
\]
The next observation explicitly relates the isotropic vectors of these forms.

\begin{obs}\label{obs:binary_iso}
Let $a, b, s, t\in \un$ and assume that $a+b$ and $as^2 + bt^2$ are nonzero. 
\begin{itemize}
\item $(x, y)$ is an isotropic vector for $\form{as^2, bt^2}$ if and only if $(sx, ty)$ is an isotropic vector for $\form{a, b}$. 
\item $(x, y)$ is an isotropic vector for $\form{a + b, (a + b)\cdot ab}$ if and only $(x - by, x + ay)$ is an isotropic vector for $\form{a, b}$.
\item $(x, y)$ is an isotropic vector for $\form{as^2 + bt^2, (as^2 + bt^2)\cdot abs^2t^2}$ if and only if $(sx - bst^2y, tx - as^2ty)$ is an isotropic vector for $\form{a, b}$.
\end{itemize}
\end{obs}

The first assertion is trivial. The second one follows by a direct calculation. The last one is an immediate consequence of the previous two. For convenience, we present one more observation that will be subsequently used in Algorithms~\ref{alg:dim5} and~\ref{alg:dim6up}.

\begin{obs}\label{obs:sign_eq}
Let $a, b\in \un$ be two nonzero elements and $\gr$ be a real place of~$K$. Set
\[
s := \bigl\lceil \abs{\sfrac{b}{a}} \bigr\rceil + 1.
\]
Then $\sgn_\gr(as^2 + b) = \sgn_\gr a$.
\end{obs}

First, we shall focus on forms of dimension~$5$. Let $\q= \form{a_1, \dotsc, a_5}$ be an isotropic form. Without loss of generality, we can assume that no $4$-dimensional subform of~$\q$ obtained by omitting one~$a_i$ is isotropic. Split $\q$ into two subforms (see also Remark~\ref{rem:permutations} below): 
\[
\q_1 := \form{a_1, a_2, a_3}
\qquad\text{and}\qquad
\q_2 := \form{a_4, a_5}.
\]
Our goal is to find an element $c\in \un$ represented by~$\q_2$ and such that $\q_1\perp \form{c}$ is isotropic. To this end, we will use Observation~\ref{obs:binary_iso} together with the Weak Approximation Theorem. The key ingredient is the following proposition.

\begin{prop}\label{prop:dim5:approximation}
For every prime $\gp\in \GP(\q_1)$ let $s_\gp, t_\gp\in \un$ be such that the form $\form{a_1, a_2, a_3, a_4s_\gp^2 + a_5t_\gp^2}$ is locally isotropic at~$\gp$. Denote the set of all these real places of~$K$ where $\q_2$ is anisotropic by $S_\infty$. Assume that for every $\gr\in S_\infty$ we have $s_\gr, t_\gr\in \un$ such that $\sgn_\gr( a_4s_\gr^2 + a_5t_\gr^2)\neq \sgn_\gr a_3$. Let $s, t\in \un$ be elements of~$K$ satisfying the following set of conditions:
\begin{enumerate}
\item\label{it:finite_prime_approx} $s\equiv s_\gp\pmod{\gp^v}$ and $t\equiv t_\gp\pmod{\gp^v}$ for all $\gp\in \GP(\q_1)$, where $v = 1 + \ord_\gp (4a_4s_\gp^2 + 4a_5 t_\gp^2)$,
\item\label{it:real_place_approx} $\abs{s - s_\gr}, \abs{t - t_\gr} < \varepsilon_\gr$ for all $\gr\in S_\infty$, where
\[
\varepsilon_\gr :=
\min\Bigl\{
\frac12, \frac{\abs{a_4s_\gr^2 + a_5t_\gr^2}}{2\cdot \abs{a_4}\cdot (1 + 2\abs{s_\gr})}, \frac{\abs{a_4s_\gr^2 + a_5t_\gr^2}}{2\cdot \abs{a_5}\cdot (1 + 2\abs{t_\gr})}
\Bigr\}.
\]
\end{enumerate}
Then the form $\form{a_1, a_2, a_3, a_4s^2 + a_5t^2}$ is isotropic.
\end{prop}

\begin{proof}
Denote $c := a_4s^2 + a_5t^2$ and $\q_* := \q_1\perp \form{c}$. We shall check that $\q_*$ is locally isotropic everywhere. It is trivially isotropic at all the complex places of~$K$. Further, if $\gr$ is a real place not in~$S_\infty$, then already $\q_1\otimes K_\gr$ is isotropic and more so $\q_*\otimes K_\gr$. Now, fix a place $\gr\in S_\infty$. We claim that the signs (with respect to~$\gr$) of $c$ and $c_\gr := a_4s_\gr^2 + a_5t_\gr^2$ coincide. Indeed, we have
\begin{align*}
c 
&= a_4\cdot (s - s_\gr + s_\gr)^2 + a_5\cdot (t - t_\gr + t_\gr)^2
\\
&= c_\gr + a_4\cdot (s - s_\gr)^2 + a_5\cdot (t - t_\gr)^2 + 2a_4\cdot (s - s_\gr)\cdot s_\gr + 2a_5\cdot (t - t_\gr)\cdot t_\gr.
\end{align*}
Observe that $\abs{s - s_\gr}^2 < \abs{s - s_\gr} < \varepsilon_\gr$, since $\varepsilon_\gr < 1$. Therefore, we have 
\begin{multline*}
\abs{ a_4\cdot (s - s_\gr)^2 + 2a_4\cdot (s - s_\gr)\cdot s_\gr }
< \abs{a_4}\cdot \varepsilon_\gr + 2\abs{a_4}\cdot \varepsilon_\gr\cdot \abs{s_\gr}
\\
\leq \abs{a_4}\cdot \frac{\abs{c_\gr}}{2\abs{a_4}\cdot (1 + 2\abs{s_\gr})}\cdot (1 + 2\abs{s_\gr})
= \frac{\abs{c_\gr}}{2}.
\end{multline*}
Analogously, we show that 
\[
\abs{a_5\cdot (t - t_\gr)^2 + 2a_5\cdot (t - t_\gr)\cdot t_\gr} < \frac{\abs{c_\gr}}{2}.
\]
It follows that $\abs{c - c_\gr} < \abs{c_\gr}$, which implies that $c$ and $c_\gr$ have equal signs, proving our claim. Now, the sign of $c_\gr$ differs from that of $a_3$. Consequently $\bigl(\q_1\perp \form{c_\gr}\bigr)\otimes K_\gr$ is isotropic and so is $\bigl(\q_1\perp \form{c}\bigr) \otimes K_\gr$.

We shall now turn our attention to non-archimedean places. If $\gp$ is a finite prime not in $\GP(\q_1)$, then already $\q_1\otimes K_\gp$ is isotropic being a $3$-dimensional form (see, e.g., \cite[Corollary~VI.2.5]{Lam2005}). Hence, $\bigl(\q_1\perp \form{c}\bigr)\otimes \Kp$ is also isotropic. Thus, the only primes we must consider are those in $\GP(\q_1)$. Fix some $\gp\in \GP(\q_1)$ and denote $c_\gp := a_4s_\gp^2 + a_5t_\gp^2$. From the congruences $s\equiv s_\gp\pmod{\gp^v}$ and $t \equiv t_\gp\pmod{\gp^v}$ we infer that 
\[
c = a_4s^2 + a_5t^2 
\equiv c_\gp
\pmod{\gp^v}.
\]
It follows from the Local Square Theorem (see, e.g., \cite[Theorem~VI.2.19]{Lam2005}) that the local square classes $c\squares{K_\gp}$ and $c_\gp\squares{K_\gp}$ coincide. Therefore, the forms $\bigl(\q_1\perp \form{c}\bigr)\otimes K_\gp$ and $\bigl(\q_1\perp \form{c_\gp}\bigr)\otimes K_\gp$ are isomorphic. In particular, $\q_1\perp \form{c}$ is locally isotropic at~$\gp$, as desired.

Now, it follows from the Strong Hasse Principle (see, e.g., \cite[Theorem VI.3.1]{Lam2005}) that $\q_1\perp \form{c}$ is isotropic over~$K$.
\end{proof}

The existence of the elements~$s$ and~$t$ in the previous proposition follows from the Weak Approximation Theorem. We can now write an explicit algorithm for finding isotropic vectors of $5$-dimensional forms.

\begin{alg}\label{alg:dim5}
Given an isotropic quadratic form $\q = \form{a_1, \dotsc, a_5}$ of dimension~$5$ over a global field~$K$, this algorithm outputs a vector $v\in K^5$ such that $q(v) = 0$. 
\begin{enumerate}
\item\label{st:dim5:quick_exit} If the $4$-dimensional form~$\q_i$, obtained from~$\q$ by omitting~$a_i$ for some $i\leq 5$, is isotropic, then execute Algorithm~\ref{alg:dim4} to construct its isotropic vector~$w$. Expand~$w$ to a $5$-dimensional vector~$v$, by inserting~$0$ at $i$-th coordinate. Output~$v$ and quit. 
\item Denote $\q_1 := \form{a_1, a_2, a_3}$. 
\item Let $S_\infty$ be the set \textup(possibly empty\textup) of all the real places~$\gr$ of~$K$ at which $\q_1$ is locally anisotropic \textup(i.e., $a_1$, $a_2$ and $a_3$ have the same signs in the completion~$K_\gr$\textup). 
\item\label{st:dim5:loop} For every $\gr\in S_\infty$:
  \begin{enumerate}
  \item If $\sgn_\gr a_4\neq \sgn_\gr a_3$, then find \textup(see Observation~\ref{obs:sign_eq}\textup) a positive integer~$s_\gr$ such that $\sgn_\gr(a_4s_\gr^2 + a_5) = \sgn_\gr a_4$ and set $t_\gr := 1$. 
  \item Otherwise, set $s_\gr := 1$ and find a positive integer $t_\gr$ such that $\sgn_\gr(a_4 + a_5t_\gr^2) = \sgn_\gr a_5$. 
  \item Set
  \[
  \varepsilon_\gr :=
  \min\Bigl\{
  \frac12, \frac{\abs{a_4s_\gr^2 + a_5t_\gr^2}}{2\cdot \abs{a_4}\cdot (1 + 2\abs{s_\gr})}, \frac{\abs{a_4s_\gr^2 + a_5t_\gr^2}}{2\cdot \abs{a_5}\cdot (1 + 2\abs{t_\gr})}
  \Bigr\}.
  \]
  \end{enumerate}
\item\label{st:dim5:local_st} For every prime $\gp\in \GP(\q_1)$ find \textup(see the comments following the algorithm\textup) a pair of elements $s_\gp, t_\gp\in \un$ such that
\[
\q_1\perp \form{ a_4s_\gp^2 + a_5t_\gp^2 }
\]
is locally isotropic at~$\gp$.
\item\label{st:dim5:wat} Using weak approximation find $s,t\in \un$ satisfying conditions \eqref{it:finite_prime_approx} and \eqref{it:real_place_approx} of Proposition~\ref{prop:dim5:approximation}.
\item\label{st:dim5:execute4} Denote $\q_* := \q_1\perp \form{a_4s^2 + a_5t^2}$. Using Algorithm~\ref{alg:dim4} find an isotropic vector $w = (w_1, w_2, w_3, w_4)$ for~$\q_*$
\item\label{st:dim5:output} Output $v := (w_1, w_2, w_3, sw_4, tw_5)$.
\end{enumerate}
\end{alg}

\begin{poc}
The correctness of the output of step~\eqref{st:dim5:quick_exit} is obvious. The Weak Approximation Theorem asserts that elements~$s$ and~$t$ in step~\eqref{st:dim5:wat} can be found. Thus, the form~$\q_*$ is isotropic by Proposition~\ref{prop:dim5:approximation}. Observation~\ref{obs:binary_iso} shows that the vector~$v$ outputted in step~\eqref{st:dim5:output} is an isotropic vector for~$\q$.
\end{poc}

It is possible to find elements $s_\gp, t_\gp\in \un$ in step~\eqref{st:dim5:local_st} of the above algorithm in a purely deterministic way. In Appendix~\ref{sec:deterministic_st}, we present one possible method to do it. However, a randomized approach is much easier to implement, and it offers a superior performance in practical tests. Unfortunately, a rigorous probabilistic analysis seems to be very difficult.

The general idea is completely straightforward. Just keep picking elements of~$\un$  at random until revealing a pair $(s_\gp, t_\gp)$ such that $(\q_1\perp \form{a_4 s_\gp^2 + a_5 t_\gp^2})\otimes \Kp$ is isotropic. If $\q_1\otimes \Kp$ is anisotropic, then by \cite[Corollary~VI.2.15]{Lam2005} it represents all local square-classes in $\sqgd[\Kp]$ except $-a_1a_2a_3\cdot \sq[\Kp]$. The group $\sqgd[\Kp]$ consists of $2^n = 4\cdot \bigl| K(\gp) \bigr|^{\ord_\gp 2}$ elements (see, e.g., \cite[Theorem~VI.2.22]{Lam2005}). In particular, if $\gp$ is non-dyadic we have $|\sqgd[\Kp]| = 4$. Therefore, if the function that maps a pair $(s, t) \in \un\times\un$ to the local square-class $(as^2 + bt^2)\cdot \sq[\Kp]$ had a uniform distribution (for some random generator of~$s$ and~$t$), the probability of failure would be $2^{-n}$. Unfortunately, it is not always the case. The reason is two-fold. For one, the form $\form{a, b}$ does not have to represent all the local square classes. It is possible to construct examples where it is the probability of \emph{success} that is as low as $2^{-n}$ (it will still be nonzero since $\q$ is isotropic by assumption) provided that all local square classes are represented by~$s$ and~$t$ with the same probability. However, here comes the second obstacle, the above assumption of equidistribution of local square classes is hard for a random generator to fulfill. Just to give the reader a glimpse of the problem, let us look at probably the simplest example when $K = \QQ$ and $\gp$ is the place associated with some odd prime number~$p\in \ZZ$. The group $\sqgd[\QQ_p]$ consists of four square classes represented by $1$, $u$, $p$ and $up$, where $(\frac{u}{p}) = -1$. Fix some $N > 0$ and a rational number $s= \sfrac{m}{n}$, where $m$ and $n$ are picked at random from the set $\{1, \dotsc, N\}$ with uniform distribution (for simplicity we ignore the sign here). Then, the probability that $s$ has an odd $p$-adic valuation is bounded from above by $2 P_1\cdot (1 - P_1)$, with an equality holding for $N$ divisible by~$p$, where
\[
P_1 = \frac{1}{p + 1}\cdot \biggl( 1 - \Bigl(-\frac{1}{p}\Bigr)^k\biggr)
\qquad\text{and}\qquad
k = \lfloor \log_p N\rfloor.
\]
Hence, in general, the two square classes of odd valuation (i.e., $p\cdot \sq[\QQ_p]$ and $up\cdot \sq[\QQ_p]$) appear far less frequently than the two others. This phenomenon is even more noticeable in general number fields. A construction of a random generator of elements of a global field uniformly distributed over local square classes exceeds the scope of this paper. A quick and dirty solution that works reasonably well in practice is to multiply with probability $\sfrac12$ elements of the field by a pre-computed uniformizer of the place in question.

It is known that algorithms relying on weak approximation can easily suffer from undesirable expression swell. This phenomenon is observed even if the approximation is limited just to finite primes (see, e.g., \cite[\S1.3.3]{Cohen2000}). It is even more pronounced when we simultaneously approximate with respect to Archimedean and non-archimedean places as we do in Algorithm~\ref{alg:dim5}. This fact stipulates the use of randomized algorithms. This approach is used in Algorithm~\ref{alg:dim5_rnd} presented below. The idea is to try to find elements~$s$ and $t$ at random instead of stitching them from local pieces using weak approximation.

\begin{alg}\label{alg:dim5_rnd}
Given an isotropic quadratic form $\q = \form{a_1, \dotsc, a_5}$ of dimension~$5$ over a global field~$K$,  and the maximal allowable number of iterations~$N$, this algorithm outputs a vector $v\in K^5$ such that $q(v) = 0$ or reports a failure. 
\begin{enumerate}
\item\label{st:dim5rnd:quick_exit} If the $4$-dimensional form~$\q_i$, obtained from~$\q$ by omitting~$a_i$ for some $i\leq 5$, is isotropic, then execute Algorithm~\ref{alg:dim4} to construct its isotropic vector~$w$. Expand~$w$ to a $5$-dimensional vector~$v$, by inserting~$0$ at $i$-th coordinate. Output~$v$ and quit. 
\item Denote $\q_1 := \form{a_1, a_2, a_3}$. 
\item\label{st:dim5rnd:loop} For $j\in \{1, \dotsc, N\}$ do:
  \begin{enumerate}
  \item Pick two random elements $s, t\in \sq$.
  \item Denote $\q_* := \q_1\perp \form{a_4s^2 + a_5t^2}$.
  \item If $\q_*$ is degenerate \textup(i.e. $a_4s^2 + a_5t^2 = 0$\textup) or anisotropic, then reiterate the loop.
  \item\label{st:dim5rnd:execute4} Using Algorithm~\ref{alg:dim4} find an isotropic vector $w = (w_1, w_2, w_3, w_4)$ for~$\q_*$
  \item Output $v := (w_1, w_2, w_3, sw_4, tw_5)$ and quit.
  \end{enumerate}
\item Report a failure.
\end{enumerate}
\end{alg}

The correctness of the above algorithm is evident. The average number of tries (hence the average running time) seems very hard to estimate for the reasons highlighted earlier. As a partial remedy in Appendix~\ref{sec:implementation}, we report the running times of (an implementation of) the algorithms presented in this paper.

We can finally turn our attention to forms of dimension~$6$ and higher. Algorithm~\ref{alg:dim6up} presented below is very similar to Algorithm~\ref{alg:dim5}. The only difference is that now, we may completely ignore finite primes. The form~$\q_*$ constructed by the algorithm has dimension $d - 1\geq 5$. Therefore, it is trivially isotropic at all the non-archimedean completions of~$K$. In particular, the next algorithm is relevant only for formally real number fields. For non-real fields (either number fields or global function fields), it suffices to take any $5$-dimensional subform of~$\q$. It will automatically be isotropic.

\begin{alg}\label{alg:dim6up}
Let $\q = \form{a_1, \dotsc, a_d}$ be an isotropic form of dimension $d\geq 6$ over a global field~$K$. This algorithm outputs a vector $v\in K^d$ such that $q(v) = 0$. 
\begin{enumerate}
\item\label{st:dim6up:quick_exit} If the $(d-1)$-dimensional form~$\q_i$, obtained from~$\q$ by omitting~$a_i$ for some index $i\leq d$, is isotropic, then construct its isotropic vector~$w$. Expand~$w$ to a $d$-dimensional vector~$v$, by inserting~$0$ at $i$-th coordinate. Output~$v$ and quit. 
\item Denote $\q_1 := \form{a_1, \dotsc, a_{d-2}}$. 
\item Let $S_\infty$ the set of all the real places~$\gr$ of~$K$, at which $\q_1$ is locally anisotropic \textup(i.e., $a_1, \dotsc, a_{d-2}$ have the same sign in the completion~$K_\gr$\textup).
\item For every $\gr\in S_\infty$:
  \begin{enumerate}
  \item If $\sgn_\gr a_{d-2}\neq \sgn_\gr a_{d-1}$, then find \textup(see Observation~\ref{obs:sign_eq}\,\textup) a positive integer~$s_\gr$ such that $\sgn_\gr(a_{d-1}s_\gr^2 + a_d) = \sgn_\gr a_{d-1}$ and set $t_\gr := 1$. 
  \item Otherwise, set $s_\gr := 1$ and find a positive integer $t_\gr$ satisfying the condition $\sgn_\gr(a_{d-1} + a_dt_\gr^2) = \sgn_\gr a_d$. 
  \item Set
  \[
  \varepsilon_\gr :=
  \min\Bigl\{
  \frac12, \frac{\abs{a_{d-1}s_\gr^2 + a_dt_\gr^2}}{2\cdot \abs{a_{d-1}}\cdot (1 + 2\abs{s_\gr})}, \frac{\abs{a_{d-1}s_\gr^2 + a_dt_\gr^2}}{2\cdot \abs{a_d}\cdot (1 + 2\abs{t_\gr})}
  \Bigr\}.
  \]
  \end{enumerate}
\item\label{st:dim6up:WAT} Using weak approximation find some $s, t\in \un$ such that 
\[
\abs{s - s_\gr} < \varepsilon_\gr
\qquad\text{and}\qquad
\abs{t - t_\gr} < \varepsilon_\gr
\]
for all real places~$\gr$ of~$K$.
\item Set $\q_* := \form{a_1, \dotsc, a_{d-2}, a_{d-1}s^2 + a_dt^2}$. Find an isotropic vector $w = (w_1, \dotsc, w_{d-1})$ of~$\q_*$.
\item Output $v = (w_1, \dotsc, w_{d-2}, sw_{d-1}, tw_{d-1})$.
\end{enumerate}
\end{alg}

\begin{rem}
Observe that the elements~$s$ and $t$ constructed in step~\eqref{st:dim5:wat} of Algorithm~\ref{alg:dim5} and step~\eqref{st:dim6up:WAT} of Algorithm~\ref{alg:dim6up} are nonzero, since the corresponding $s_\gr$'s and $t_\gr$'s are nonzero integers and $\abs{s - s_\gr}, \abs{t - t_\gr} < \sfrac12$. This fact guarantees that the form~$\q_*$ is non-degenerate.
\end{rem}

\begin{rem}\label{rem:permutations}
Algorithms \ref{alg:dim5}--\ref{alg:dim6up} rely on finding an element~$c$ represented by a binary subform of~$\q$ and such that $\q_1\perp \form{c}$ is isotropic. All three algorithms construct the binary form from the last two coefficients of~$\q$. The only reason for such an arbitrary partitioning of the list of coefficients is the clarity of presentation. Obviously, we can take any binary subform $\form{a_i, a_j}$ of~$\q$ instead. The form~$\q_1$ is then given by the remaining $(d - 2)$ coefficients. It opens a way for some possible optimizations that may reduce the size of resulting elements~$s$ and~$t$. One natural strategy is to select the indices $i$ and $j$ in a way that minimizes the number of real places at which $\q_1$ is anisotropic.
\end{rem}

\subsection*{Acknowledgments} The author wishes to thank Radan Kučera for helpful comments on an early version of this paper and the anonymous reviewer for suggesting a simpler version of algorithms for forms of dimension~$5$ and above.

\bibliographystyle{plain} 
\bibliography{isovec}

\begin{thebibliography}{10}

\bibitem{BPR2003}
Saugata Basu, Richard Pollack, and Marie-Fran\c{c}oise Roy.
\newblock {\em Algorithms in real algebraic geometry}, volume~10 of {\em
  Algorithms and Computation in Mathematics}.
\newblock Springer-Verlag, Berlin, 2003.

\bibitem{Belabas2004}
Karim Belabas.
\newblock Topics in computational algebraic number theory.
\newblock {\em J. Th\'{e}or. Nombres Bordeaux}, 16(1):19--63, 2004.

\bibitem{BF2014}
Jean-Fran\c{c}ois Biasse and Claus Fieker.
\newblock Subexponential class group and unit group computation in large degree
  number fields.
\newblock {\em LMS J. Comput. Math.}, 17(suppl. A):385--403, 2014.

\bibitem{BRS2015}
Alex Bocharov, Martin Roetteler, and Krysta~M. Svore.
\newblock Efficient synthesis of probabilistic quantum circuits with fallback.
\newblock {\em Phys. Rev. A}, 91:052317, May 2015.

\bibitem{BCP1997}
Wieb Bosma, John Cannon, and Catherine Playoust.
\newblock The {M}agma algebra system. {I}. {T}he user language.
\newblock {\em J. Symbolic Comput.}, 24(3-4):235--265, 1997.
\newblock Computational algebra and number theory (London, 1993).

\bibitem{BLP1993}
Joe~P. Buhler, Hendrik~W. Lenstra, Jr., and Carl Pomerance.
\newblock Factoring integers with the number field sieve.
\newblock In {\em The development of the number field sieve}, volume 1554 of
  {\em Lecture Notes in Math.}, pages 50--94. Springer, Berlin, 1993.

\bibitem{CBFS2023}
John Cannon, Wieb Bosma, Claus Fieker, and Allan~Steel (eds.).
\newblock {\em Handbook of Magma Functions}, 2.28 edition, 2023.

\bibitem{Castel2013}
Pierre Castel.
\newblock Solving quadratic equations in dimension 5 or more without factoring.
\newblock In {\em A{NTS} {X}---{P}roceedings of the {T}enth {A}lgorithmic
  {N}umber {T}heory {S}ymposium}, volume~1 of {\em Open Book Ser.}, pages
  213--233. Math. Sci. Publ., Berkeley, CA, 2013.
\newblock \url{http://dx.doi.org/10.2140/obs.2013.1.213}.

\bibitem{CM1998}
Todd Cochrane and Patrick Mitchell.
\newblock Small solutions of the {L}egendre equation.
\newblock {\em J. Number Theory}, 70(1):62--66, 1998.

\bibitem{Cohen1993}
Henri Cohen.
\newblock {\em A course in computational algebraic number theory}, volume 138
  of {\em Graduate Texts in Mathematics}.
\newblock Springer-Verlag, Berlin, 1993.
\newblock \url{https://doi.org/10.1007/978-3-662-02945-9}.

\bibitem{Cohen2000}
Henri Cohen.
\newblock {\em Advanced topics in computational number theory}, volume 193 of
  {\em Graduate Texts in Mathematics}.
\newblock Springer-Verlag, New York, 2000.
\newblock \url{https://doi.org/10.1007/978-1-4419-8489-0}.

\bibitem{CR2003}
John Cremona and David Rusin.
\newblock Efficient solution of rational conics.
\newblock {\em Math. Comp.}, 72(243):1417--1441 (electronic), 2003.
\newblock \url{http://dx.doi.org/10.1090/S0025-5718-02-01480-1}.

\bibitem{DMR2021}
Mawunyo~Kofi Darkey-Mensah and Beata Rothkegel.
\newblock Computing the length of sum of squares and {P}ythagoras element in a
  global field.
\newblock {\em Fund. Inform.}, 184(4):297--306, 2021.

\bibitem{FJP1997}
Claus Fieker, Andreas Jurk, and Michael~E. Pohst.
\newblock On solving relative norm equations in algebraic number fields.
\newblock {\em Math. Comp.}, 66(217):399--410, 1997.
\newblock \url{https://doi.org/10.1090/S0025-5718-97-00761-8}.

\bibitem{FP1983}
Ulrich Fincke and Michael~E. Pohst.
\newblock A procedure for determining algebraic integers of given norm.
\newblock In {\em Computer algebra ({L}ondon, 1983)}, volume 162 of {\em
  Lecture Notes in Comput. Sci.}, pages 194--202. Springer, Berlin, 1983.
\newblock \url{https://doi.org/10.1007/3-540-12868-9_103}.

\bibitem{Garbanati1980}
Dennis~A. Garbanati.
\newblock An algorithm for finding an algebraic number whose norm is a given
  rational number.
\newblock {\em J. Reine Angew. Math.}, 316:1--13, 1980.
\newblock \url{https://doi.org/10.1515/crll.1980.316.1}.

\bibitem{GMN2013}
Jordi Gu\`ardia, Jes\'{u}s Montes, and Enric Nart.
\newblock A new computational approach to ideal theory in number fields.
\newblock {\em Found. Comput. Math.}, 13(5):729--762, 2013.
\newblock \url{https://doi.org/10.1007/s10208-012-9137-5}.

\bibitem{IKR2019}
G\'{a}bor Ivanyos, P\'{e}ter Kutas, and Lajos R\'{o}nyai.
\newblock Explicit equivalence of quadratic forms over {$\mathbb{F}_q(t)$}.
\newblock {\em Finite Fields Appl.}, 55:33--63, 2019.
\newblock \url{https://doi.org/10.1016/j.ffa.2018.09.003}.

\bibitem{UP2024}
Jorge Jim\'{e}nez~Urroz and Jacek Pomyka\l~a.
\newblock Factoring numbers with elliptic curves.
\newblock {\em Ramanujan J.}, 64(1):265--273, 2024.

\bibitem{KNW2019}
Habiba Kadiri, Nathan Ng, and Peng-Jie Wong.
\newblock The least prime ideal in the {C}hebotarev density theorem.
\newblock {\em Proc. Amer. Math. Soc.}, 147(6):2289--2303, 2019.

\bibitem{Koprowski2021}
Przemysław Koprowski.
\newblock Computing singular elements modulo squares.
\newblock {\em Fund. Inform.}, 179(3):227--238, 2021.

\bibitem{KR2023}
Przemysław Koprowski and Beata Rothkegel.
\newblock The anisotropic part of a quadratic form over a number field.
\newblock {\em J. Symbolic Comput.}, 115:39--52, 2023.

\bibitem{KMZC2022}
P\'{e}ter Kutas, Micka\"{e}l Montessinos, Gergely Z\'{a}br\'{a}di, and
  T\'{\i}mea Csah\'{o}k.
\newblock Finding nontrivial zeros of quadratic forms over rational function
  fields of characteristic 2.
\newblock In {\em I{SSAC} '22---{P}roceedings of the 2022 {I}nternational
  {S}ymposium on {S}ymbolic and {A}lgebraic {C}omputation}, pages 235--244.
  ACM, New York, [2022] \copyright 2022.

\bibitem{LMO1979}
Jeffrey~C. Lagarias, Hugh~L. Montgomery, and Andrew~M. Odlyzko.
\newblock A bound for the least prime ideal in the {C}hebotarev density
  theorem.
\newblock {\em Invent. Math.}, 54(3):271--296, 1979.

\bibitem{Lam2005}
Tsit-Yuen Lam.
\newblock {\em Introduction to quadratic forms over fields}, volume~67 of {\em
  Graduate Studies in Mathematics}.
\newblock American Mathematical Society, Providence, RI, 2005.

\bibitem{LV2018}
Jonathan~D. Lee and Ramarathnam Venkatesan.
\newblock Rigorous analysis of a randomised number field sieve.
\newblock {\em J. Number Theory}, 187:92--159, 2018.

\bibitem{LJS2024}
Sumin Leem, Michael~J. Jacobson, and Renate Scheidler.
\newblock Solving norm equations in global function fields.
\newblock {\em Res. Number Theory}, 11(1):Paper No. 17, 2025.

\bibitem{LW1992}
David~B. Leep and Adrian~R. Wadsworth.
\newblock The {H}asse norm theorem mod squares.
\newblock {\em J. Number Theory}, 42(3):337--348, 1992.
\newblock \url{https://doi.org/10.1016/0022-314X(92)90098-A}.

\bibitem{Lenstra1987}
Hendrik~W. Lenstra, Jr.
\newblock Factoring integers with elliptic curves.
\newblock {\em Ann. of Math. (2)}, 126(3):649--673, 1987.

\bibitem{lmfdb}
The {LMFDB Collaboration}.
\newblock The {L}-functions and modular forms database.
\newblock \url{https://www.lmfdb.org}, 2025.
\newblock [Online; accessed 10 March 2025].

\bibitem{MKV2023}
Victor Magron, Przemysław Koprowski, and Tristan Vaccon.
\newblock Pourchet's theorem in action: decomposing univariate nonnegative
  polynomials as sums of five squares.
\newblock In {\em Proceedings of the {I}nternational {S}ymposium on {S}ymbolic
  \& {A}lgebraic {C}omputation ({ISSAC} 2023)}, pages 425--433. ACM, New York,
  [2023] \copyright 2023.

\bibitem{OMeara2000}
O.~Timothy O'Meara.
\newblock {\em Introduction to quadratic forms}.
\newblock Classics in Mathematics. Springer-Verlag, Berlin, 2000.
\newblock Reprint of the 1973 edition.

\bibitem{PTBW2020}
Lillian~B. Pierce, Caroline~L. Turnage-Butterbaugh, and Melanie~Matchett Wood.
\newblock An effective {C}hebotarev density theorem for families of number
  fields, with an application to {$\ell$}-torsion in class groups.
\newblock {\em Invent. Math.}, 219(2):701--778, 2020.

\bibitem{Pohst2000}
Michael~E. Pohst.
\newblock On {L}egendre's equation over number fields.
\newblock volume~56, pages 535--546. 2000.
\newblock Dedicated to Professor K\'{a}lm\'{a}n Gy\H{o}ry on the occasion of
  his 60th birthday.

\bibitem{Schicho1998}
Josef Schicho.
\newblock Rational parameterization of real algebraic surfaces.
\newblock In {\em Proceedings of the 1998 {I}nternational {S}ymposium on
  {S}ymbolic and {A}lgebraic {C}omputation ({R}ostock)}, pages 302--308. ACM,
  New York, 1998.

\bibitem{Schicho2000}
Josef Schicho.
\newblock Proper parametrization of surfaces with a rational pencil.
\newblock In {\em Proceedings of the 2000 {I}nternational {S}ymposium on
  {S}ymbolic and {A}lgebraic {C}omputation ({S}t. {A}ndrews)}, pages 292--300.
  ACM, New York, 2000.

\bibitem{Simon2002}
Denis Simon.
\newblock Solving norm equations in relative number fields using {$S$}-units.
\newblock {\em Math. Comp.}, 71(239):1287--1305, 2002.
\newblock \url{https://doi.org/10.1090/S0025-5718-02-01309-1}.

\bibitem{Simon2005}
Denis Simon.
\newblock Solving quadratic equations using reduced unimodular quadratic forms.
\newblock {\em Math. Comp.}, 74(251):1531--1543 (electronic), 2005.
\newblock \url{http://dx.doi.org/10.1090/S0025-5718-05-01729-1}.

\bibitem{Smart1998}
Nigel~P. Smart.
\newblock {\em The algorithmic resolution of {D}iophantine equations},
  volume~41 of {\em London Mathematical Society Student Texts}.
\newblock Cambridge University Press, Cambridge, 1998.

\bibitem{Szymiczek1997}
Kazimierz Szymiczek.
\newblock {\em Bilinear algebra}, volume~7 of {\em Algebra, Logic and
  Applications}.
\newblock Gordon and Breach Science Publishers, Amsterdam, 1997.
\newblock An introduction to the algebraic theory of quadratic forms.

\bibitem{vHC2006}
Mark van Hoeij and John Cremona.
\newblock Solving conics over function fields.
\newblock {\em J. Th\'{e}or. Nombres Bordeaux}, 18(3):595--606, 2006.
\newblock \url{http://jtnb.cedram.org/item?id=JTNB_2006__18_3_595_0}.

\bibitem{Voight2013}
John Voight.
\newblock Identifying the matrix ring: algorithms for quaternion algebras and
  quadratic forms.
\newblock In {\em Quadratic and higher degree forms}, volume~31 of {\em Dev.
  Math.}, pages 255--298. Springer, New York, 2013.

\bibitem{GG2013}
Joachim von~zur Gathen and J\"{u}rgen Gerhard.
\newblock {\em Modern computer algebra}.
\newblock Cambridge University Press, Cambridge, third edition, 2013.

\end{thebibliography}

\appendix
\section{Complexity of Algorithm~\ref{alg:dim4}}\label{sec:complexity}
A detailed complexity analysis of algorithms presented in this paper seems to be very difficult, and definitely exceeds the competence of the author. One of the reasons is that to the best of our knowledge not all of the external subroutines, we rely on, offer such a rigorous complexity analysis. In this appendix we provide a crude estimate for the time complexity of Algorithm~\ref{alg:dim4}.

Let $K$ be a global field represented either in the form $K = \QQ(\theta)$, if $K$ is a number field, or $K = \FF_{p^k}(x, \theta)$, if it is a function field. Denote the monic minimal polynomial of~$\theta$ by~$f$ and let $d$ be the degree of~$f$. An element of~$K$ can be expressed as a linear combination of successive powers of~$\theta$, say $a = a_0 + a_1\theta + \dotsb + a_{d-1}\theta^{d-1}$. Define the size of~$a$ to be the sum of bit-sizes of $a_0, \dotsc, a_{d-1}$. For a quadratic form~$\q$ over~$K$, define its size as a sum of sizes of its coefficients.

In order to construct the set~$\GP(\q)$ in step~\eqref{st:dim4:initilization} of the algorithm, one has to factor the principal ideals generated by $K$-integral representatives of square classes of coefficients of~$\q$. If $K$ is a number field, this problem can be reduced in polynomial time to the problem of factoring integers (see, e.g., \cite[\S2.3.5]{Cohen2000}). The best known algorithms for the latter task are sub-exponential in size of the input. GNFS (see, e.g., \cite{BLP1993}, \cite{LV2018}) has complexity $L_n(\sfrac13, (\sfrac{64}{9})^{1/3})$, while ECM (see, e.g., \cite{Lenstra1987}, \cite{UP2024}) has complexity $L_p(\sfrac12, \sqrt{2})$, where $p$ is the smallest prime dividing the integer to be factored. Now, the computation of a norm increases the size of the input by the factor $|K:\QQ|$ and so the complexity of step~\eqref{st:dim4:initilization} for number fields is sub-exponential in the size of the form and the degree of the field.

If $K$ is a global function field, a finite extension of $\FF_{p^k}(x)$, then the problem of factoring the ideals can be analogously reduced to factoring polynomials over a finite field. The later problem is well-known to be polynomial-time in the size of the input (see, e.g., \cite[Section~14.8]{GG2013}) Hence, for function fields step~\eqref{st:dim4:initilization} has polynomial time in size of the input and the degree~$d$.

Now, steps (\ref{st:dim4:real_places}--\ref{st:dim4:vec_u}) are meaningful only for formally real fields (hence, necessarily number fields). They boil down to performing a real root isolation of polynomials with rational coefficients. This can be done in polynomial time with respect to the degree $|K:\QQ|$ and the logarithm of the root distance. For details we refer the reader to \cite[Section~10.2]{BPR2003} (see also \cite[Section~3]{Belabas2004}).

As expected, it is the analysis of the loop in step~\eqref{st:dim4:loop} which is most involving. The construction in step~\eqref{st:dim4:basis} of the basis~$\SB$ from scratch during the first execution of the loop can be reduced in polynomial time (with respect to the size of~$S$) to the problem of computing the $S$-unit group. For a general number field, the author is not aware of any algorithm for a latter task that would offer a complexity better than exponential in the size of~$S$ and the bit-size of the discriminant of~$K$. Related sub-exponential algorithms (see, e.g., \cite{BF2014}) rely on some heuristics. On the other hand, as explained in Remark~\ref{rem:incremental}, in subsequent iterations of the loop, after expanding~$S$, the basis~$\SB$ can be updated incrementally in time polynomial with respect to the size of~$S$ and the class number of~$K$. The corresponding algorithm for function fields is described in \cite{Koprowski2021}. It is likely that one can use a similar incremental algorithm also for number fields. However, to the best of our knowledge no such algorithm has ever been published.

Computations of Hilbert symbols in steps (\ref{st:dim4:vec_v}, \ref{st:dim4:vec_w}) and (\ref{st:dim4:mtx_B}, \ref{st:dim4:mtx_C}) can be done using Voight's algorithm (see \cite{Voight2013}) which works in the polynomial time. The previous comments on the complexity of steps (\ref{st:dim4:real_places}--\ref{st:dim4:vec_u}) apply also to step~\eqref{st:dim4:mtx_A}. All the remaining steps of the loop use only linear algebra over~$\FF_2$, and so they run in the polynomial time with respect to the size of $S\cup S_\infty$. 

It remains to estimate the needed number of iterations of loop~\eqref{st:dim4:loop}. Denote by~$\GQ$ the set of all these primes~$\gq$ for which system~\eqref{eq:quaternary} has a solution in $\Sing{\GP(q)\cup\{\gq\}}$. As explained in Remark~\ref{rem:density} after Algorithm~\ref{alg:dim4}, the set~$\GQ$ has a positive density. In fact, if one repeats all the steps of the proof of \cite[Lemma~2.1]{LW1992}, keeping track of the index of the used congruence subgroup, one can actually show that this density is bounded from below by
\[
\prod_{\gp\in \GP(\q)\cup \{\gp_0\}} \hspace*{-2mm}(N\gp)^{-1 - \ord_\gp 4},
\]
where $\gp_0$ is some fixed (arbitrarily chosen) prime not in~$\GP(\q)$. 

As mentioned in Remark~\ref{rem:exhaustive_search}, Algorithm~\ref{alg:dim4} can be viewed either as a randomized or a deterministic one. In the former case, the average number of iterations of loop~\eqref{st:dim4:loop} is bounded by the inverse of the above product. Consequently, it is polynomial in the degree of the field and the size of coefficients. Summarizing, if Algorithm~\ref{alg:dim4} is treated as a randomized one, then loop~\eqref{st:dim4:loop} has an exponential time but only due to the first execution of step~\eqref{st:dim4:basis}. If the first construction of the basis~$\SB$ is moved to the initialization phase before the loop (which would simplify the complexity analysis, but obscure the presentation of the algorithm), then step~\eqref{st:dim4:loop} would execute in a randomized polynomial time.

Alternatively, we may treat Algorithm~\ref{alg:dim4} as a deterministic one, exhaustively iterating over the primes of~$K$. Although it may be tempting to base our analysis in this setting on an effective version of the Chebotarev Density Theorem (see, e.g., \cite{PTBW2020} and the citations wherein), it is not an optimal strategy. The reason is that we are not interested in the number of primes in~$\GQ$, whose norms sit in a given interval. What we need is to estimate how quickly we will find the first prime (with respect to the order in which we scan them) sitting in~$\GQ$. This is sometimes referred to as the ``first prime problem''. It is known to work in polynomial time with respect to the discriminant of~$K$ by \cite{LMO1979} (see also \cite{KNW2019} and citations within). Consequently, step~\eqref{st:dim4:loop} runs in the deterministic polynomial time with respect to the the degree of~$K$, its discriminant and size of the coefficients, provided that the first construction of the basis~$\SB$ is moved to the initialization phase (see the comments of the end of the previous paragraph).

Finally, let us look at step~\eqref{st:dim4:solve_ternary}. Here, we need to find isotropic vectors of two ternary forms. As explained in Section~\ref{sec:dim3}, this problem is equivalent to solving two norm equations. Algorithms devoted to this task has been recalled in the introduction. For number fields, Simon's algorithm \cite{Simon2002} seems to be the most popular one. (It can be also adapted to function fields.) The author is not aware of any complete rigorous complexity analysis of Simon's algorithm. It is claimed (without proof) in \cite{BRS2015} that over cyclotomic fields a solution of the norm equation can be computed in probabilistic polynomial time with respect to the bit-size of the input, provided that the factorization of the input is known. The latter condition can be safely assumed in our case since the set~$S$ is given and we are interested only in solutions modulo squares. Yet still, this claim refers only to cyclotomic fields. Conversely, in \cite{MKV2023} there is a complete complexity analysis of a stripped down version of Simon's algorithm, showing that the computation time is polynomial with respect to the size of the representation for the system of fundamental units of the unit group of the upper field, which in turn is exponential with respect to the bit-size of the coefficients (and also the degree of the extension, but in our case it is always~$2$). On the other hand, for global function fields, Leem, Jacobson and Scheidler (see \cite{LJS2024}) presented recently a new algorithm (in fact two new algorithms) for solving norm equations. Their algorithm works in polynomial time with respect to the size of the underlying finite field and the bit-size of coefficients, and in sub-exponential time with respect to the genus of the field. (It is also exponential with respect to the degree of the field extension but again it is always $2$ in our case.)

\section{Experimental results}\label{sec:implementation}
The algorithms presented in this paper were implemented by the author as part of the CQF library for the computer algebra system Magma \cite{BCP1997}. The purpose of these implementations was not only to provide a proof of concept but also to evaluate their runtime performance on real-world hardware.

Figures~\ref{fig:dim4_bitsize}--\ref{fig:dim4-5_char} present the measured execution times for constructing isotropic vectors of quadratic forms of dimensions~$4$ and $5$, analyzed with respect to the factors identified in Appendix~\ref{sec:complexity}: the bit-size of the form, the discriminant of the base field, the degree of the base field, and the size of the underlying finite field (in the case of global function fields). In each case, four distinct execution times are reported:
\begin{itemize}
\item The median execution time (computed over multiple samples) for the complete construction of an isotropic vector.
\item The median execution time for the initialization phase, preceding the loop \eqref{st:dim4:loop} in Algorithm~\ref{alg:dim4} (or loop~\eqref{st:dim5rnd:loop} in Algorithm~\ref{alg:dim5_rnd}). This phase primarily involves constructing the set~$\GP(\q)$, i.e., factoring the principal ideals generated by the form’s coefficients.
\item The median execution time of the main loop, i.e., loop~\eqref{st:dim4:loop} in Algorithm~\ref{alg:dim4} or loop~\eqref{st:dim5rnd:loop} in Algorithm~\ref{alg:dim5_rnd}.
\item The median execution time of solving the norm equations (step~\eqref{st:dim4:solve_ternary} in Algorithm~\ref{alg:dim4}) or of finding an isotropic vector for a $4$-dimensional form (step~\eqref{st:dim5rnd:execute4} in Algorithm~\ref{alg:dim5_rnd}).
\end{itemize}
All the measurements were done on Intel\textregistered{} Core\texttrademark~i5-9400F with 32GB of RAM running Magma~2.28-6. The CQF library, along with the Magma scripts used for testing, input data, scripts generating the input data, and detailed execution timings, are available for download from the author’s homepage: \url{http://www.pkoprowski.eu/cqf}.

All the forms used for the test were randomly generated in such a way that the quick-exit test in step~\eqref{st:dim4:quick_exit} of Algorithm~\ref{alg:dim4} (resp. Algorithm~\ref{alg:dim5_rnd}) is never triggered.

The reader may wonder why we report medians rather than the more commonly used average execution times. This choice is motivated by the fact that, among seemingly similar forms, there may be instances where the execution time differs drastically from the rest. The presence of such outliers can significantly distort the mean (since execution times cannot be negative, atypical values always skew the mean upwards).

This phenomenon is illustrated in Table~\ref{tab:dim4_deg14}, which presents the actual measured execution times from test number~4 (described in detail below) for Algorithm~\ref{alg:dim4} applied to forms over $20$ different number fields of degree~$14$. Notably, for all form-field pairs, the execution time of the loop in step~\eqref{st:dim4:loop} is approximately 2 seconds, except for the form over the field LMFDB 14.0.318497840669359.1, where this step takes nearly 90 seconds. As a result, the median execution times (used in Figure~\ref{fig:dim4-5_degree}) provide a more robust measure, minimizing the impact of extreme values. Specifically, the medians are: $40.46$, $0.30$, $2.30$, and $38.25$, while the averages are $77.01$, $0.65$, $6.76$, and $69.59$, respectively. The downside of using the medians is that the medians of the partial times do not sum up to the median of the entire execution times.
\begin{table}
\caption{\label{tab:dim4_deg14}Detailed execution times of Algorithm~\ref{alg:dim4} for forms over number fields   of degree~$14$.}
\begin{tabular}{cr@{.}lr@{.}lr@{.}lr@{.}l}
LMFDB field label & 
\multicolumn{2}{c}{\begin{tabular}{@{}c@{}}total \\ time\end{tabular}} &  
\multicolumn{2}{c}{\begin{tabular}{@{}c@{}}initia-\\ lization\end{tabular}} &
\multicolumn{2}{c}{loop~\eqref{st:dim4:loop}} &
\multicolumn{2}{c}{step~\eqref{st:dim4:solve_ternary}}\\
\hline\rule{0mm}{4mm}
14.2.116768679686993.1 &   39 & 93 & \rule{3mm}{0mm}0 & 13 &  1 & 60 &  38 & 20\\
14.0.118236562059083.1 &  314 & 23 & 0 & 34 &  2 & 16 & 311 & 73\\
14.4.297261105546875.1 &   33 & 20 & 0 & 12 &  1 & 88 &  31 & 20\\
14.0.280155320935227.1 &   42 & 62 & 0 & 27 &  4 & 05 &  38 & 30\\
14.2.103886084489093.1 &   31 & 72 & 0 & 17 &  2 & 29 &  55 & 00\\
14.0.318497840669359.1 &   48 & 97 & 0 & 30 &  3 & 20 &  45 & 47\\
14.0.263657868936003.1 &   34 & 89 & 0 & 14 &  1 & 74 &  33 & 01\\
14.2.212764932023296.1 &   34 & 13 & 0 & 30 &  3 & 15 &  30 & 68\\
14.2.245928712415581.1 &   31 & 56 & 5 & 07 &  2 & 38 &  24 & 11\\
14.2.245928712415581.1 &   39 & 89 & 0 & 58 &  3 & 54 &  35 & 77\\
14.0.172575877556547.1 &   41 & 58 & 0 & 09 &  2 & 29 &  39 & 20\\
14.0.118236562059083.1 &   34 & 24 & 0 & 93 &  2 & 46 &  30 & 85\\
14.0.280155320935227.1 &   40 & 99 & 0 & 34 &  2 & 32 &  38 & 33\\
14.0.263657868936003.1 &   75 & 04 & 0 & 08 &  1 & 63 &  73 & 33\\
14.4.297261105546875.1 &   54 & 04 & 0 & 14 &  1 & 75 &  52 & 15\\
14.4.334624313614643.1 &   29 & 38 & 0 & 30 &  3 & 07 &  26 & 01\\
14.0.318497840669359.1 &  384 & 87 & 0 & 38 & 89 & 92 & 294 & 57\\
14.2.154291569457229.1 &   27 & 51 & 0 & 92 &  1 & 76 &  24 & 83\\
14.4.297261105546875.1 &  142 & 08 & 0 & 17 &  1 & 86 & 140 & 05\\
14.0.318497840669359.1 &   59 & 38 & 2 & 32 &  2 & 50 &  54 & 56\\
\end{tabular} 
\end{table} 

In the first test, we measured how the execution time of Algorithm~\ref{alg:dim4} depends on the bit-size of the input form. To this end, we considered successive intervals of length 64 bits: $[1, 64], [65, 128], \dotsc, [960, 1024]$. For each interval, we generated $64$ random forms of dimension~$4$ and with bit-sizes in the given interval. Then, we measured the four execution times specified above. For the base fields we used two number fields: $\QQ(i)$ and LMFDB~8.4.15908237.1, and two function fields: $\FF_7(x)$ and $\FF_{11}[x,y]/(y^2 - x^3 - 3x + 1)$. An example form (of size $867$ bits) used in the test is: 
\begin{align*}
\langle 
& 
\frac{- 227\,737\,278\,468\,653\,418\,985\,585\,661\,180\,683\,861\,540\,509\,597}
{156\,827\,860\,436\,626\,092\,968\,994\,140\,625\,000} 
\\&\hspace{20mm}+ 
\frac{952\,381\,959\,201\,357\,284\,651\,722\,474\,420\,956\,860\,826\,197\,846\, i }
{156\,827\,860\,436\,626\,092\,968\,994\,140\,625\,000}, 
\\&\hspace*{10mm}
- 1\,210 - 1\,100\, i,
\\&\hspace*{10mm}
-\frac{16\,425\,718\,989\,170\,711\,280\,553 - 16\,163\,790\,306\,359\,826\,077\,304\,i}
{39\,601\,497\,501\,562\,500},
\\&\hspace*{10mm}
-\frac{57\,431\,391\,141 + 140\,480\,602\,822\, i}{199\,001\,250} 
\rangle
\end{align*}
Figure~\ref{fig:dim4_bitsize} presents the medians over successive intervals of the measured timings. For the above example form the initialization phase took $4.06$\;s, the execution of loop~\eqref{st:dim4:loop} took $0.08$\;s, and the solution of the two norm equations in step~\eqref{st:dim4:solve_ternary} needed $0.62$\;s. This sums up to the total execution time of $4.76$\;s.

Observe that in Figure~\eqref{subfig:GFF_bitsize}, for forms over $\FF_{11}[x,y]/(y^2 - x^3 - 3x + 1)$, we present the execution times of step~\eqref{st:dim4:solve_ternary}, and consequently also the full execution times, only for the first two intervals. Our implementation relies on Magma's built-in method for finding rational points on cubic curves (equivalent to solving norm equations). For forms over $128$ bits with coefficients in this function field, the execution time of step~\eqref{st:dim4:solve_ternary} tends to exceed the threshold of half an hour for a single form, that we set to finish the tests in a reasonable time frame. It is planned to replace Magma's built-in method for solving norm equations over function fields with our own implementation in a future version of CQF.

The second test is very similar to the first one, except here we measured the execution time of Algorithm~\ref{alg:dim5_rnd}. Hence, for the same intervals that we used in the previous test, we generated quintic isotropic forms. This time, we used $32$ forms for each interval. As explained above, the forms used in the test are constructed in a way that avoids the quick-exit in step~\eqref{st:dim5rnd:quick_exit}. However, such forms are rather rare. In particular, we failed to find even a single such form of size below $64$ bits over \emph{any} of the fields used in the test. Actually the minimal bit-sizes for $5$-dimensional isotropic forms we found, for which the test in step~\eqref{st:dim5rnd:quick_exit} is not triggered, are: $100$ bits for $\QQ(i)$, $219$ bits for LMFDB~8.4.15908237.1, $130$ bits for $\FF_7(x)$ and $160$ bits for $\FF_{11}[x,y]/(y^2 - x^3 - 3x + 1)$. Figure~\ref{fig:dim5_bitsize} shows the median execution times of this test. Again, in the last (sub)figure, we omit the execution time of step~\eqref{st:dim5rnd:execute4} for the same reason as before. The reader should notice that for $5$-dimensional forms, the total running time is completely dominated by the construction of an isotropic vector of a $4$-dimensional form in step~\eqref{st:dim5rnd:execute4}.

In the next test, we checked how the discriminant of the ground field impacts the execution times of Algorithms~\ref{alg:dim4} and~\ref{alg:dim5_rnd}. To this end, from the popular number-theoretic library LMFDB \cite{lmfdb}, we downloaded a list of the quadratic fields with discriminants ranging from $-10^5$ to $10^5$. Then, we considered successive intervals $[1, 2\,000], [2\,001, 4\,000], \dotsc,  [98\,001, 100\,000]$. For each interval, we constructed $50$ random pairs, consisting of a field~$K$ with $|\disc K|$ in the specified interval and an isotropic quadratic form over~$K$. The test was then repeated four times. Twice for forms of dimension~$4$ and sizes in the intervals $[112, 144]$ bits and $[240, 272]$ bits, respectively. Next, twice for forms of dimension~$5$ and bit-sizes in the same intervals. The resulting execution times are presented in Figure~\ref{fig:dim4-5_disc}. Once again, for $5$-dimensional forms, the total execution time is completely dominated by the construction of an isotropic vector of the $4$-dimensional form.

The fourth test shows how the execution times of the two algorithms depend on the degree of the base field. We again used the LMFDB database to obtain a list of number fields of degrees up to~$15$ and with the root discriminants in the interval $[10, 11]$. Then, for each degree, we constructed $20$ random pairs consisting of a field with the specified degree and a quartic isotropic form over this field with size between $512$ and $768$ bits (resp. $512$ and $1\,024$ bits for quintic forms). The median execution times depending on the degree are presented in Figure~\ref{subfig:dim4_degree}. The test was then repeated for $5$- dimensional forms, except that here we performed the tests only for degrees up to~$12$. The results are presented in Figure~\ref{subfig:deg5_degree}. It should be noticed that with the increasing degree, it is the solution of the associated norm equations in step~\eqref{st:dim4:solve_ternary}, which becomes the dominating step.

In the last test, we analyze how the size of the underlying finite field impacts the execution times of the two algorithms. For the reason explained in the first test, we restricted ourselves only to forms over rational function fields. For every prime number~$p$ ranging from~$3$ to~$337$, we constructed the rational function fields $\FF_p(x)$ and $50$ random isotropic forms of dimension~$4$ (respectively $5$) over it. The bit-sizes of the forms were restricted to the interval $[ 490, 544]$ bits. Figure~\ref{fig:dim4-5_char} shows the running times of Algorithms~\ref{alg:dim4} and~\ref{alg:dim5_rnd} (and their respective parts). Observe that for $4$-dimensional forms, it is loop~\ref{st:dim4:loop} that dominates the execution time, while for $5$- dimensional forms, almost the whole time is spent on finding an isotopic vector of the resulting $4$-dimensional form.

\tikzset{
  FullStyle/.style={lightgray, line width = 2bp, line join=round},
  PqStyle/.style={densely dashed, line join=round},
  LoopStyle/.style={solid, line join=round},
  NeqStyle/.style={dotted, line width = 1bp, line join=round}
}
%
%
%
\begin{figure}[p]
\begin{subfigure}{.5\textwidth}
  \centering
  \begin{tikzpicture}
  \begin{axis}[
     width=65mm,
      xlabel={size of the form [bits]},
      ylabel={time [s]},
      ylabel near ticks,
      xmin = 0, xmax = 1024,
      xtick distance = 256,
      ymax = 5.5,
  ]
  \addplot [FullStyle] table [
    x expr= \thisrow{bs}-64, y expr=\thisrow{med1}
  ] {results_dim4_bitsize_Qi.dat};
  \addplot [PqStyle] table [
    x expr= \thisrow{bs}-64, y expr=\thisrow{med2}
  ] {results_dim4_bitsize_Qi.dat};
  \addplot [LoopStyle] table [
    x expr= \thisrow{bs}-64, y expr=\thisrow{med3}
  ] {results_dim4_bitsize_Qi.dat};
  \addplot [NeqStyle] table [
    x expr= \thisrow{bs}-64, y expr=\thisrow{med4}
  ] {results_dim4_bitsize_Qi.dat};
  \end{axis}
  \end{tikzpicture}
  \caption{Base field $\QQ(i)$}
\end{subfigure}%
\begin{subfigure}{.5\textwidth}
  \centering
  \begin{tikzpicture}
  \begin{axis}[
     width=65mm,
      xlabel={size of the form [bits]},
      grid style=dashed,
      xmin = 0, xmax = 1024,
      xtick distance = 256,
      ymax = 5.5,
  ]
  \addplot [FullStyle] table [
    x expr= \thisrow{bs}, y expr=\thisrow{med1}
  ] {results_dim4_bitsize_LMDB_8.4.15908237.1.dat};
  \addplot [PqStyle] table [
    x expr= \thisrow{bs}, y expr=\thisrow{med2}
  ] {results_dim4_bitsize_LMDB_8.4.15908237.1.dat};
  \addplot [LoopStyle] table [
    x expr= \thisrow{bs}, y expr=\thisrow{med3}
  ] {results_dim4_bitsize_LMDB_8.4.15908237.1.dat};
  \addplot [NeqStyle] table [
    x expr= \thisrow{bs}, y expr=\thisrow{med4}
  ] {results_dim4_bitsize_LMDB_8.4.15908237.1.dat};
  \end{axis}
  \end{tikzpicture}
  \caption{Base field LMDB 8.4.15908237.1}
\end{subfigure}
\bigskip

\begin{subfigure}{.5\textwidth}
  \centering
  \begin{tikzpicture}
  \begin{axis}[
     width=65mm,
      xlabel={size of the form [bits]},
      ylabel={time [s]},
      ylabel near ticks,
      xmin = 0, xmax = 1024,
      xtick distance = 256,
      ymax = 3.5,
  ]
  \addplot [FullStyle] table [
    x expr= \thisrow{bs}-64, y expr=\thisrow{med1}
  ] {results_dim4_RFF_bitsize_FF7.dat};
  \addplot [PqStyle] table [
    x expr= \thisrow{bs}-64, y expr=\thisrow{med2}
  ] {results_dim4_RFF_bitsize_FF7.dat};
  \addplot [LoopStyle] table [
    x expr= \thisrow{bs}-64, y expr=\thisrow{med3}
  ] {results_dim4_RFF_bitsize_FF7.dat};
  \addplot [NeqStyle] table [
    x expr= \thisrow{bs}-64, y expr=\thisrow{med4}
  ] {results_dim4_RFF_bitsize_FF7.dat};
  \end{axis}
  \end{tikzpicture}
  \caption{Base field $\FF_7(x)$}
\end{subfigure}%
\begin{subfigure}{.5\textwidth}
  \centering
  \begin{tikzpicture}
  \begin{axis}[
     width=65mm,
      xlabel={size of the form [bits]},
      grid style=dashed,
      xmin = 0, xmax = 1024,
      xtick distance = 256,
      ymax = 3.5,
  ]
  \addplot [FullStyle] table [
    x expr= \thisrow{bs}, y expr=\thisrow{med1}
  ] {results_dim4_GFF_bitsize.dat};
  \addplot [PqStyle] table [
    x expr= \thisrow{bs}, y expr=\thisrow{med2}
  ] {results_dim4_GFF_bitsize.dat};
  \addplot [LoopStyle] table [
    x expr= \thisrow{bs}, y expr=\thisrow{med3}
  ] {results_dim4_GFF_bitsize.dat};
  \addplot [NeqStyle] table [
    x expr= \thisrow{bs}, y expr=\thisrow{med4}
  ] {results_dim4_GFF_bitsize.dat};
  \end{axis}
  \end{tikzpicture}
  \caption{\label{subfig:GFF_bitsize}Base field $\FF_{11}[x,y]/(y^2 - x^3 - 3x + 1)$}
\end{subfigure}
\caption{Median running times of Algorithm~\ref{alg:dim4} depending on the bit-size of the form. The thick gray line represents the total execution time. The dashed line is the time of the initialization phase. The thin solid line is the time spent in loop~\eqref{st:dim4:loop}. Finally, the dotted line shows the time used to solve norm equations in step~\eqref{st:dim4:solve_ternary}.}
\label{fig:dim4_bitsize}
\end{figure}

%
%
%
\begin{figure}[p]
\begin{subfigure}{.5\textwidth}
  \centering
  \begin{tikzpicture}
  \begin{axis}[
     width=65mm,
      xlabel={size of the form [bits]},
      ylabel={time [s]},
      ylabel near ticks,
      xmin = 0, xmax = 1024,
      xtick distance = 256,
      ymax = 35,
  ]
  \addplot [FullStyle] table [
    x expr= \thisrow{bs}-64, y expr=\thisrow{med1}
  ] {results_dim5_bitsize_Qi.dat};
  \addplot [PqStyle] table [
    x expr= \thisrow{bs}-64, y expr=\thisrow{med2}
  ] {results_dim5_bitsize_Qi.dat};
  \addplot [LoopStyle] table [
    x expr= \thisrow{bs}-64, y expr=\thisrow{med3}
  ] {results_dim5_bitsize_Qi.dat};
  \addplot [NeqStyle] table [
    x expr= \thisrow{bs}-64, y expr=\thisrow{med4}
  ] {results_dim5_bitsize_Qi.dat};
  \end{axis}
  \end{tikzpicture}
  \caption{Base field $\QQ(i)$}
\end{subfigure}%
\begin{subfigure}{.5\textwidth}
  \centering
  \begin{tikzpicture}
  \begin{axis}[
     width=65mm,
      xlabel={size of the form [bits]},
      grid style=dashed,
      xmin = 0, xmax = 1024,
      xtick distance = 256,
      ymax = 35,
  ]
  \addplot [FullStyle] table [
    x expr= \thisrow{bs} - 64, y expr=\thisrow{med1}
  ] {results_dim5_bitsize_LMFDB_8.4.15908237.1.dat};
  \addplot [PqStyle] table [
    x expr= \thisrow{bs} - 64, y expr=\thisrow{med2}
  ] {results_dim5_bitsize_LMFDB_8.4.15908237.1.dat};
  \addplot [LoopStyle] table [
    x expr= \thisrow{bs} - 64, y expr=\thisrow{med3}
  ] {results_dim5_bitsize_LMFDB_8.4.15908237.1.dat};
  \addplot [NeqStyle] table [
    x expr= \thisrow{bs} - 64, y expr=\thisrow{med4}
  ] {results_dim5_bitsize_LMFDB_8.4.15908237.1.dat};
  \end{axis}
  \end{tikzpicture}
  \caption{Base field LMDB 8.4.15908237.1}
\end{subfigure}
\bigskip

\begin{subfigure}{.5\textwidth}
  \centering
  \begin{tikzpicture}
  \begin{axis}[
     width=65mm,
      xlabel={size of the form [bits]},
      ylabel={time [s]},
      ylabel near ticks,
      xmin = 0, xmax = 1024,
      xtick distance = 256,
      ymax = 35,
  ]
  \addplot [FullStyle] table [
    x expr= \thisrow{bs} - 64, y expr=\thisrow{med1}
  ] {results_dim5_RFF_bitsize.dat};
  \addplot [PqStyle] table [
    x expr= \thisrow{bs} - 64, y expr=\thisrow{med2}
  ] {results_dim5_RFF_bitsize.dat};
  \addplot [LoopStyle] table [
    x expr= \thisrow{bs} - 64, y expr=\thisrow{med3}
  ] {results_dim5_RFF_bitsize.dat};
  \addplot [NeqStyle] table [
    x expr= \thisrow{bs} - 64, y expr=\thisrow{med4}
  ] {results_dim5_RFF_bitsize.dat};
  \end{axis}
  \end{tikzpicture}
  \caption{Base field $\FF_7(x)$}
\end{subfigure}%
\begin{subfigure}{.5\textwidth}
  \centering
  \begin{tikzpicture}
  \begin{axis}[
     width=65mm,
      xlabel={size of the form [bits]},
      grid style=dashed,
      xmin = 0, xmax = 1024,
      xtick distance = 256,
      ymax = 1.1,
  ]
  \addplot [PqStyle] table [
    x expr= \thisrow{bs}, y expr=\thisrow{med2}
  ] {results_dim5_GFF_bitsize.dat};
  \addplot [LoopStyle] table [
    x expr= \thisrow{bs}, y expr=\thisrow{med3}
  ] {results_dim5_GFF_bitsize.dat};
  \end{axis}
  \end{tikzpicture}
  \caption{Base field $\FF_{11}[x,y]/(y^2 - x^3 - 3x + 1)$}
\end{subfigure}
\caption{Median running times of Algorithm~\ref{alg:dim5_rnd} depending on the bit-size of the form. The thick gray line represents the total execution time. The dashed line is the time of the initialization phase. The thin solid line is the time spent in loop~\eqref{st:dim5rnd:loop}. Finally, the dotted line shows the time used to find an isotropic vector of the $4$-dimensional form in step~\eqref{st:dim5rnd:execute4}.}
\label{fig:dim5_bitsize}
\end{figure}

%
%
\begin{figure}[p]
\begin{subfigure}{.5\textwidth}
  \centering
  \begin{tikzpicture}
  \begin{axis}[
     width=65mm,
      xlabel={discriminant},
      ylabel={time [s]},
      ylabel near ticks,
      ymax = 1.85,
  ]
  \addplot [FullStyle] table [
    x expr= \thisrow{disc}, y expr=\thisrow{med1}
  ] {results_dim4_disc_112-144bits.dat};
  \addplot [PqStyle] table [
    x expr= \thisrow{disc}, y expr=\thisrow{med2}
  ] {results_dim4_disc_112-144bits.dat};
  \addplot [LoopStyle] table [
    x expr= \thisrow{disc}, y expr=\thisrow{med3}
  ] {results_dim4_disc_112-144bits.dat};
  \addplot [NeqStyle] table [
    x expr= \thisrow{disc}, y expr=\thisrow{med4}
  ] {results_dim4_disc_112-144bits.dat};
  \end{axis}
  \end{tikzpicture}
  \caption{$\dim \q = 4$, size $[112, 144]$ bits}
\end{subfigure}%
\begin{subfigure}{.5\textwidth}
  \centering
  \begin{tikzpicture}
  \begin{axis}[
     width=65mm,
      xlabel={discriminant},
      grid style=dashed,
      ymax = 1.85,
  ]
  \addplot [FullStyle] table [
    x expr= \thisrow{disc}, y expr=\thisrow{med1}
  ] {results_dim4_disc_240-272bits.dat};
  \addplot [PqStyle] table [
    x expr= \thisrow{disc}, y expr=\thisrow{med2}
  ] {results_dim4_disc_240-272bits.dat};
  \addplot [LoopStyle] table [
    x expr= \thisrow{disc}, y expr=\thisrow{med3}
  ] {results_dim4_disc_240-272bits.dat};
  \addplot [NeqStyle] table [
    x expr= \thisrow{disc}, y expr=\thisrow{med4}
  ] {results_dim4_disc_240-272bits.dat};
  \end{axis}
  \end{tikzpicture}
  \caption{$\dim \q = 4$, size $[240, 272]$ bits}
\end{subfigure}
\bigskip

\begin{subfigure}{.5\textwidth}
  \centering
  \begin{tikzpicture}
  \begin{axis}[
     width=65mm,
      xlabel={discriminant},
      ylabel={time [s]},
      ylabel near ticks,
      ymax = 6.5,
  ]
  \addplot [FullStyle] table [
    x expr= \thisrow{disc}, y expr=\thisrow{med1}
  ] {results_dim5_disc_112-144bits.dat};
  \addplot [PqStyle] table [
    x expr= \thisrow{disc}, y expr=\thisrow{med2}
  ] {results_dim5_disc_112-144bits.dat};
  \addplot [LoopStyle] table [
    x expr= \thisrow{disc}, y expr=\thisrow{med3}
  ] {results_dim5_disc_112-144bits.dat};
  \addplot [NeqStyle] table [
    x expr= \thisrow{disc}, y expr=\thisrow{med4}
  ] {results_dim5_disc_112-144bits.dat};
  \end{axis}
  \end{tikzpicture}
  \caption{$\dim \q = 5$, size $[112, 144]$ bits}
\end{subfigure}%
\begin{subfigure}{.5\textwidth}
  \centering
  \begin{tikzpicture}
  \begin{axis}[
     width=65mm,
      xlabel={discriminant},
      grid style=dashed,
      ymax = 6.5,
  ]
  \addplot [FullStyle] table [
    x expr= \thisrow{disc}, y expr=\thisrow{med1}
  ] {results_dim5_disc_240-272bits.dat};
  \addplot [PqStyle] table [
    x expr= \thisrow{disc}, y expr=\thisrow{med2}
  ] {results_dim5_disc_240-272bits.dat};
  \addplot [LoopStyle] table [
    x expr= \thisrow{disc}, y expr=\thisrow{med3}
  ] {results_dim5_disc_240-272bits.dat};
  \addplot [NeqStyle] table [
    x expr= \thisrow{disc}, y expr=\thisrow{med4}
  ] {results_dim5_disc_240-272bits.dat};
  \end{axis}
  \end{tikzpicture}
  \caption{$\dim \q = 5$, size $[240, 272]$ bits}
\end{subfigure}
\caption{Median running times of Algorithms~\ref{alg:dim4} (top row) and~\ref{alg:dim5_rnd} (bottom row) for number fields depending on the discriminant of the field. The thick gray line represents the total execution time. The dashed line is the time of the initialization phase. The thin solid line is the time spent in loop~\eqref{st:dim4:loop} of Algorithm~\ref{alg:dim4} (top), or loop~\eqref{st:dim5rnd:loop} of Algorithm~\ref{alg:dim5_rnd} (bottom). Finally, the dotted line shows the time used to find an isotropic vector of a lower-dimensional form.}
\label{fig:dim4-5_disc}
\end{figure}

%
%
\begin{figure}[p]
\begin{subfigure}{.5\textwidth}
  \centering
  \begin{tikzpicture}
  \begin{axis}[
     width=65mm,
      xlabel={degree $|K:\QQ|$},
      ylabel={time [s]},
      ylabel near ticks,
      xmin = 1, xmax = 16,
      ymax = 85,
  ]
  \addplot [FullStyle] table [
    x expr= \thisrow{deg}, y expr=\thisrow{med1}
  ] {results_dim4_deg_rd_10-11.dat};
  \addplot [PqStyle] table [
    x expr= \thisrow{deg}, y expr=\thisrow{med2}
  ] {results_dim4_deg_rd_10-11.dat};
  \addplot [LoopStyle] table [
    x expr= \thisrow{deg}, y expr=\thisrow{med3}
  ] {results_dim4_deg_rd_10-11.dat};
  \addplot [NeqStyle] table [
    x expr= \thisrow{deg}, y expr=\thisrow{med4}
  ] {results_dim4_deg_rd_10-11.dat};
  \end{axis}
  \end{tikzpicture}
  \caption{\label{subfig:dim4_degree}$\dim\q = 4$}
\end{subfigure}%
\begin{subfigure}{.5\textwidth}
  \centering
  \begin{tikzpicture}
  \begin{axis}[
      width=65mm,
      xlabel={degree $|K:\QQ|$},
      xmin = 1, xmax = 13,
      xtick distance = 4,
      ymax = 85,
  ]
  \addplot [FullStyle] table [
    x expr= \thisrow{deg}, y expr=\thisrow{med1}
  ] {results_dim5_deg_rd_10-11.dat};
  \addplot [PqStyle] table [
    x expr= \thisrow{deg}, y expr=\thisrow{med2}
  ] {results_dim5_deg_rd_10-11.dat};
  \addplot [LoopStyle] table [
    x expr= \thisrow{deg}, y expr=\thisrow{med3}
  ] {results_dim5_deg_rd_10-11.dat};
  \addplot [NeqStyle] table [
    x expr= \thisrow{deg}, y expr=\thisrow{med4}
  ] {results_dim5_deg_rd_10-11.dat};
  \end{axis}
  \end{tikzpicture}
  \caption{\label{subfig:deg5_degree}$\dim\q = 5$}
\end{subfigure}
\caption{Median running times of Algorithms~\ref{alg:dim4} (left) and Algorithm~\ref{alg:dim5_rnd} (right) depending on the degree of the base field. The thick gray line represents the total execution time. The dashed line is the time of the initialization phase. The thin solid line is the time spent in  loop~\eqref{st:dim4:loop} of Algorithm~\ref{alg:dim4} (left), or loop~\eqref{st:dim5rnd:loop} of Algorithm~\ref{alg:dim5_rnd} (right). Finally, the dotted line shows the time used to find an isotropic vector of the lower-dimensional form.}
\label{fig:dim4-5_degree}
\end{figure}

%
%
\begin{figure}[p]
\begin{subfigure}{.5\textwidth}
  \centering
  \begin{tikzpicture}
  \begin{axis}[
     width=65mm,
      xlabel={$p$},
      ylabel={time [s]},
      ylabel near ticks,
      ymax = 4.5,
  ]
  \addplot [FullStyle] table [
    x expr= \thisrow{p}, y expr=\thisrow{med1}
  ] {results_dim4_RFF_basefield_480-544bits.dat};
  \addplot [PqStyle] table [
    x expr= \thisrow{p}, y expr=\thisrow{med2}
  ] {results_dim4_RFF_basefield_480-544bits.dat};
  \addplot [LoopStyle] table [
    x expr= \thisrow{p}, y expr=\thisrow{med3}
  ] {results_dim4_RFF_basefield_480-544bits.dat};
  \addplot [NeqStyle] table [
    x expr= \thisrow{p}, y expr=\thisrow{med4}
  ] {results_dim4_RFF_basefield_480-544bits.dat};
  \end{axis}
  \end{tikzpicture}
  \caption{$\dim\q = 4$}
\end{subfigure}%
\begin{subfigure}{.5\textwidth}
  \centering
  \begin{tikzpicture}
  \begin{axis}[
     width=65mm,
      xlabel={$p$},
      grid style=dashed,
      ymax = 4.5,
  ]
  \addplot [FullStyle] table [
    x expr= \thisrow{p}, y expr=\thisrow{med1}
  ] {results_dim5_RFF_basefield_480-544bits.dat};
  \addplot [PqStyle] table [
    x expr= \thisrow{p}, y expr=\thisrow{med2}
  ] {results_dim5_RFF_basefield_480-544bits.dat};
  \addplot [LoopStyle] table [
    x expr= \thisrow{p}, y expr=\thisrow{med3}
  ] {results_dim5_RFF_basefield_480-544bits.dat};
  \addplot [NeqStyle] table [
    x expr= \thisrow{p}, y expr=\thisrow{med4}
  ] {results_dim5_RFF_basefield_480-544bits.dat};
  \end{axis}
  \end{tikzpicture}
  \caption{$\dim \q = 5$}
\end{subfigure}
\caption{Median running times of Algorithms~\ref{alg:dim4} (left) and \ref{alg:dim5_rnd} (right) for rational function fields over finite fields depending on the size of the field of constants. The thick gray line represents the total execution time. The dashed line is the time of the initialization phase. The thin solid line is the time spent in loop~\eqref{st:dim4:loop} of Algorithm~\ref{alg:dim4} (left), or loop~\eqref{st:dim5rnd:loop} of Algorithm~\ref{alg:dim5_rnd} (right). Finally, the dotted line shows the time used to find an isotropic vector of the lower-dimensional form.}
\label{fig:dim4-5_char}
\end{figure}

\section{Deterministic construction of elements $s$ and $t$}\label{sec:deterministic_st}
As explained in the comments following Algorithm~\ref{alg:dim5}, an efficient and recommended method to construct elements $s_\gp, t_\gp\in \un$ needed in step~\eqref{st:dim5:local_st} of Algorithm~\ref{alg:dim5} is by a random search. Nonetheless, in the author's opinion, it is in good taste to show that we are not limited to stochastic arguments but can construct such elements also using a purely deterministic method. Below, we will treat dyadic and non-dyadic primes uniformly here. However, if one is interested exclusively in the non-dyadic case (e.g., if one works over global function fields), then there are some possible optimizations that are not covered here. We begin with an auxiliary algorithm.

\begin{alg}\label{alg:st_units}
Let $\gp$ be a prime of a number field~$K$. Given three nonzero elements $a$, $b$ and $d$ of~$K$ of valuations zero, this algorithm outputs $s, t\in \un$ such that $as^2 + bt^2\not\equiv d\pmod{\sq[\Kp]}$.
\begin{enumerate}
\item\label{st:stdu_a_vs_d} If $a\cdot d$ is not a local square at~$\gp$, then:
  \begin{enumerate}
  \item Find a uniformizer~$\pi$ of~$\gp$. 
  \item Output $s := 1$, $t := 2\pi$ and quit. 
  \end{enumerate}
\item\label{st:stdu_b_vs_d} If $b\cdot d$ is not a local square at~$\gp$, then 
  \begin{enumerate}
  \item Find a uniformizer~$\pi$ of~$\gp$. 
  \item Output $s := 2\pi$, $t := 1$ and quit.
  \end{enumerate}
\item If all three elements $a$, $b$ and $d$ belong to the same local-square class, then:
  \begin{enumerate}
  \item Find $u$ such that $a = bu^2$. 
  \item Find $s_*, t_*\in \un$ such that $s_*^2 + t_*^2$ is not a square in~$\Kp$ \textup(use, e.g., \cite{DMR2021}\textup). 
  \item Output $s := s_*$ and $t := ut_*$.
  \end{enumerate}
\end{enumerate}
\end{alg}

\begin{poc}
First, suppose that $a\cdot d\notin \sq[\Kp]$. Thus, 
\[
as^2 + bt^2 
= a + b\cdot 4\pi^2
\equiv a\pmod{\gp^{1 + \ord_\gp 4}}.
\]
It follows from the Local Square Theorem that $as^2 + bt^2$ sits in the local square class of~$a$, which differs from the local square class of~$d$ since $a\cdot d$ is not a local square. It proves the correctness of step~\eqref{st:stdu_a_vs_d}. The correctness of step~\eqref{st:stdu_b_vs_d} is completely analogous. Finally, consider the case when $a$, $b$, and $d$ belong to the same local square class. Since $s_*^2 + t_*^2$ is not a square, we have
\[
as^2 + bt^2 = a\cdot (s_*^2 + t_*^2)\not\equiv a\equiv d\pmod{\sq[\Kp]}.
\qedhere
\]
\end{poc}

Now, we can present a deterministic algorithm for constructing~$s$ and~$t$.

\begin{alg}
Let $\gp$ be a prime of a number field~$K$. Given three nonzero elements $a$, $b$ and $d$ of~$K$, this algorithm outputs $s, t\in \un$ such the $as^2 + bt^2\not\equiv d \pmod{\sq[\Kp]}$.
\begin{enumerate}
\item If all three elements $a$, $b$ and $d$ have the same parity of valuation, then: 
  \begin{enumerate}
  \item Find a uniformizer~$\pi$ of~$\gp$.
  \item Write $a$, $b$ and $d$ as
  \[
  a = a_*\cdot \pi^{\ord_\gp a},\quad
  b = b_*\cdot \pi^{\ord_\gp b},\quad
  d = d_*\cdot \pi^{\ord_\gp d},
  \]
  for some $\gp$-units $a_*, b_*, d_*\in \un$.
  \item Execute Algorithm~\ref{alg:st_units} with input $a_*$, $b_*$, $d_*$, and denote its output by $s_*$ and $t_*$.
  \item\label{st:stnd:units} Output $s := s_*\cdot \pi^{\lceil \sfrac{\ord_\gp d}{2}\rceil - \lceil \sfrac{\ord_\gp a}{2}\rceil}$, $t := t_*\cdot \pi^{\lceil \sfrac{\ord_\gp d}{2}\rceil - \lceil \sfrac{\ord_\gp b}{2}\rceil}$ and quit.
  \end{enumerate}
\item\label{st:stnd:a_vs_d} If $\ord_\gp d \not\equiv \ord_\gp a \pmod{2}$, then set $s := 1$ and let $t := \pi^k$ be such that $2k > \ord_\gp a + \ord_\gp b$. 
\item Otherwise, if $\ord_\gp d \not\equiv \ord_\gp b \pmod{2}$, then set $t := 1$ and let $s := \pi^k$ be such that $2k > ord_\gp a + \ord_\gp b$. 
\item Output $s$ and $t$.
\end{enumerate}
\end{alg}

\begin{poc}
If all three elements have the same parity of valuations, then the correctness of the result outputted in step~\eqref{st:stnd:units} follows from the correctness of Algorithm~\ref{alg:st_units}. Conversely, suppose that the valuations of~$a$ and~$d$ have different parity. Then for~$s$ and~$t$ constructed in step~\eqref{st:stnd:a_vs_d} we have 
\[
\ord_\gp (as^2 + bt^2) 
= \ord_\gp as^2 
\equiv \ord_\gp a
\not\equiv \ord_\gp d\pmod{2}. 
\]
It follows that local square classes of~$d$ and $as^2 + bt^2$ differ. One proceeds analogously when the parity of $\ord_\gp b$ and $\ord_\gp d$ differ.
\end{poc}

\end{document}